\documentclass[11pt]{article}
\usepackage{fullpage,graphicx,amssymb,amsmath,psfrag,color}
\definecolor{ddarkbrown}{rgb}{0.5,0.2,0.05}
\definecolor{bbluegray}{rgb}{0.05,0,0.5}
\usepackage[round]{natbib}
\usepackage[colorlinks,citecolor=bbluegray,linkcolor=ddarkbrown,urlcolor=blue,breaklinks]{hyperref}
\usepackage{pdfsync}
\usepackage{times}
\usepackage{algorithmic,algorithm}

\newcommand{\BEAS}{\begin{eqnarray*}}
\newcommand{\EEAS}{\end{eqnarray*}}
\newcommand{\BEA}{\begin{eqnarray}}
\newcommand{\EEA}{\end{eqnarray}}
\newcommand{\BEQ}{\begin{equation}}
\newcommand{\EEQ}{\end{equation}}
\newcommand{\BIT}{\begin{itemize}}
\newcommand{\EIT}{\end{itemize}}
\newcommand{\BNUM}{\begin{enumerate}}
\newcommand{\ENUM}{\end{enumerate}}

\newcommand{\BA}{\begin{array}}
\newcommand{\EA}{\end{array}}

\newcommand{\eg}{{\it e.g.}}

\newcommand{\ones}{\mathbf 1}

\newcommand{\reals}{{\mbox{\bf R}}}

\newcommand{\symm}{{\mbox{\bf S}}}  

\newcommand{\Rank}{\mathop{\bf Rank}}

\newcommand{\Card}{\mathop{\bf Card}}
\newcommand{\Tr}{\mathop{\bf Tr}}
\newcommand{\diag}{\mathop{\bf diag}}
\newcommand{\lambdamax}{{\lambda_{\rm max}}}

\newcommand{\idm}{\mathbf{I}}

\newcommand{\QED}{~~\rule[-1pt]{6pt}{6pt}}

\newcommand{\dsp}{\displaystyle}

\newcommand{\argmax}{\mathop{\rm argmax}}

\newtheorem{theorem}{Theorem}

\newtheorem{proposition}[theorem]{Proposition}
\newtheorem{remark}[theorem]{Remark}

\newcounter{exno}

\newenvironment{proof}{\textbf{Proof.}}{\QED\bigskip}

{\begin{quote}}{\end{quote}}

\makeatletter
\long\def\@makecaption#1#2{
   \vskip 9pt 
   \begin{small}
   \setbox\@tempboxa\hbox{{\bf #1:} #2}
   \ifdim \wd\@tempboxa > 5.5in
        \begin{center}
        \begin{minipage}[t]{5.5in}
        \addtolength{\baselineskip}{-0.95pt}
        {\bf #1:} #2 \par
        \addtolength{\baselineskip}{0.95pt}
        \end{minipage}
        \end{center}
   \else 
    \hbox to\hsize{\hfil\box\@tempboxa\hfil}  
   \fi
   \end{small}\par
}
\makeatother

\newcounter{oursection}

\newcounter{lecture}

\let\oldmarginpar\marginpar
\renewcommand\marginpar[1]{\-\oldmarginpar[\raggedleft\footnotesize #1]%
{\raggedright\footnotesize #1}}

\begin{document}
\title{Sparse PCA: Convex Relaxations, Algorithms and Applications}
\author{Youwei Zhang\thanks{EECS, U.C. Berkeley. Berke.ley, CA 94720. \texttt{zyw@eecs.berkeley.edu}}, ~Alexandre d'Aspremont\thanks{ORFE, Princeton University, Princeton, NJ 08544. \texttt{aspremon@princeton.edu}}, ~Laurent El Ghaoui\thanks{EECS, U.C. Berkeley. Berkeley, CA 94720. \texttt{elghaoui@eecs.berkeley.edu}}}
\maketitle

\begin{abstract}
Given a sample covariance matrix, we examine the problem of maximizing the variance explained by a linear combination of the input variables while constraining the number of nonzero coefficients in this combination. This is known as sparse principal component analysis and has a wide array of applications in machine learning and engineering. Unfortunately, this problem is also combinatorially hard and we discuss convex relaxation techniques that efficiently produce good approximate solutions. We then describe several algorithms solving these relaxations as well as greedy algorithms that iteratively improve the solution quality. Finally, we illustrate sparse PCA in several applications, ranging from senate voting and finance to news data.
\end{abstract}

\section{Introduction}
Principal component analysis (PCA) is a classical tool for data analysis, visualization and dimensionality reduction and has a wide range of applications throughout science and engineering. Starting from a multivariate data set, PCA finds linear combinations of the variables called \emph{principal components}, corresponding to orthogonal directions maximizing variance in the data.
Numerically, a full PCA involves a singular value decomposition of the data matrix.

One of the key shortcomings of PCA is that the factors are linear combinations of \emph{all} original variables; that is, most of factor coefficients (or loadings) are non-zero.  This means that while PCA facilitates model interpretation and visualization by concentrating the information in a few factors, the factors themselves are still constructed using all variables, hence are often hard to interpret.

In many applications, the coordinate axes involved in the factors have a direct physical interpretation. In financial or  biological applications, each axis might correspond to a specific asset or gene. In problems such as these, it is natural to seek a trade-off between the two goals of \emph{statistical fidelity} (explaining most of the variance in the data) and \emph{interpretability} (making sure that the factors involve only a few coordinate axes).  Solutions that have only a few nonzero coefficients in the principal components are usually easier to interpret. Moreover, in some applications, nonzero coefficients have a direct cost (\eg, transaction costs in finance) hence there may be a direct trade-off between statistical fidelity and practicality. Our aim here is to efficiently derive \emph{sparse principal components}, i.e, a set of sparse vectors that explain a maximum amount of  variance. Our motivation is that in many applications, the decrease in statistical fidelity required to obtain sparse factors is  small and relatively benign. 

In what follows, we will focus on the problem of finding sparse factors which explain a maximum amount of variance in the original data, which can be written
\BEQ \label{eq:spca-intro}
\max_{ \| z \|  \leq 1 } z^T \Sigma z - \rho \Card(z)
\EEQ
in the variable $z\in\reals^n$, where $\Sigma\in\symm_n$ is the (symmetric positive semi-definite) sample covariance matrix, $\rho$ is a parameter controlling sparsity, and $\Card(z)$ denotes the cardinal (or $\ell_0$ norm) of $z$, i.e. the number of non zero coefficients of $z$.

While PCA is numerically easy, each factor requires computing a leading eigenvector, which can be done in $O(n^2)$ floating point operations using the Lanczos method for example (see e.g. \cite[\S8.3, \S9.1.1]{Golu90} or \cite{Saad92} for details), sparse PCA is a hard combinatorial problem. In fact, \cite{Mogh06a} show that the subset selection problem for ordinary least squares, which is NP-hard~\cite{Nata95}, can be reduced to a sparse generalized eigenvalue problem, of which sparse PCA is a particular instance. Sometimes factor rotation techniques are used to post-process the results from PCA and improve interpretability (see QUARTIMAX by \cite{Neuh54}, VARIMAX by \cite{Kais58} or \cite{Joll95} for a discussion). Results by e.g. \cite{Amin08} show that this naive approach has significantly worst convergence rates than the relaxations we present here. Another straightforward solution is to \emph{threshold} to zero loadings with small magnitude \cite{cadi95}, but outside of easy cases, the methods highlighted below always perform better in situation when only a few observations are available or when significant noise is present.

A more systematic approach to the problem arose in recent years, with various researchers proposing nonconvex algorithms  (e.g., SCoTLASS by \cite{Joll03}, SLRA by \cite{Zhan02a} or D.C. based methods \cite{Srip07} which find modified principal components  with zero loadings. The SPCA algorithm, which is based on the representation of PCA as a regression-type optimization problem~\cite{Zou04}, allows the application of the LASSO \cite{tibs96}, a penalization technique based on the $\ell_1$ norm. With the exception of simple thresholding, all the algorithms above require solving non convex problems. Recently also, \cite{dasp04a} derived an $\ell_1$ based semidefinite relaxation for the sparse PCA problem (\ref{eq:spca-intro}) with a complexity of $O(n^{4}\sqrt{\log n})$ for a given $\rho$. \cite{Mogh06b} used greedy search and branch-and-bound methods to solve small instances of problem (\ref{eq:spca-intro}) exactly and get good solutions for larger ones. Each step of this greedy algorithm has complexity $O(n^3)$, leading to a total complexity of $O(n^4)$ for a full set of solutions. \cite{Mogh07} improve this bound in the regression/discrimination case. \cite{Jour08} use an extension of the power method to (locally) solve the problem defined here, as well as the ``block'' problem of finding several sparse principal components at once. Loss of orthogonality means that there is no natural method for deflating the matrix once a sparse principal component is found and \cite{Mack09} discusses several options, comparing the variance vs. orthogonality/sparsity tradeoffs they imply. Finally, \cite{Amin08} derive explicit sample size thresholds for recovery of true sparse vector using either simple thresholding methods or semidefinite relaxations, in a spiked model for the covariance. 

Here, we detail two semidefinite relaxations for sparse PCA, and describe algorithms to solve the relaxations efficiently. We also test these techniques on various data sets: newsgroup data, Senate voting records and stock market returns.

\subsection*{Notation}
For a vector $z\in\reals$, we let $\|z\|_1 = \sum_{i=1}^n |z_i|$ and $\|z\| = \left( \sum_{i=1}^n z_i^2 \right)^{1/2}$, $\Card(z)$ is the cardinality of $z$, i.e. the number of nonzero coefficients of $z$, while the support $I$ of $z$ is the set $\{i:~z_i\neq 0\}$ and we use $I^c$ to denote its complement. For $\beta\in\reals$, we write $\beta_+ = \max \{\beta,0\}$ and for $X\in\symm_n$ (the set of symmetric matrix of size $n\times n$) with eigenvalues $\lambda_i$, $\Tr(X)_+=\sum_{i=1}^n \max\{\lambda_i,0\}$. The vector of all ones is written $\ones$, while the identity matrix is written $\idm$. The diagonal matrix with the vector $u$ on the diagonal is written $\diag(u)$.

\section{Semidefinite relaxations} \label{s:sdp-relax}
\label{sec:spca} Let $\Sigma\in\symm_n$ be a symmetric matrix. We consider the following sparse PCA problem
\BEQ \label{eq:pca-card}
\phi(\rho) \equiv \max_{ \| z \|  \leq 1 } z^T \Sigma z - \rho \Card(z)
\EEQ
in the variable $z\in\reals^n$ where $\rho >0$ is a parameter controlling sparsity. We assume without loss of generality that $\Sigma\in\symm_n$ is positive semidefinite and that the $n$ variables are ordered by decreasing marginal variances, i.e.
that $\Sigma_{11} \geq \ldots \geq \Sigma_{nn}$. We also assume that we are given a square root $A$ of the matrix $\Sigma$ with $ \Sigma = A^T A$, where $A\in\reals^{n \times n}$ and we denote by $a_1,\dots,a_n \in \reals^n$ the columns of $A$. Note that the problem and our algorithms are invariant by permutations of $\Sigma$ and by the choice of square root $A$. In practice, we are very often given the data matrix $A$ instead of the covariance $\Sigma$.

A problem that is directly related to (\ref{eq:pca-card}) is that of computing a cardinality constrained maximum eigenvalue, by solving
\BEQ\label{eq:card-eig}
\BA{ll}
\mbox{maximize} & z^T \Sigma z\\
\mbox{subject to} & \Card(z) \leq k\\
& \|z\|=1,
\EA\EEQ
in the variable $z\in\reals^n$. Of course, this problem and (\ref{eq:pca-card}) are related. By weak duality, an upper bound on the optimal value of (\ref{eq:card-eig}) is given by
\[
\inf_{\rho \in P} \phi(\rho) + \rho k.
\]
where $P$ is the set of penalty values for which $\phi(\rho)$ has been computed. This means in particular that if a point $z$ is provably optimal for (\ref{eq:pca-card}), it is also globally optimum for (\ref{eq:card-eig}) with $k=\Card(z)$.

\subsection{A semidefinite relaxation with $\ell_1$ penalization}\label{ss:l1-relax}
Here, we briefly recall the $\ell_1$ based relaxation derived in \cite{dasp04a}. Following the \emph{lifting procedure} for semidefinite relaxation described in \cite{Lova91}, \cite{Aliz95}, \cite{Lema99} for example, we rewrite (\ref{eq:card-eig}) as
\BEQ \label{eq-variat-matrix-prog}
\BA{ll}
\mbox{maximize} & \Tr(\Sigma X)\\
\mbox{subject to} & \Tr(X)=1\\
& \Card(X) \leq k^2\\
& X \succeq 0,~\Rank(X)=1,
\EA \EEQ
in the (matrix) variable $X\in \symm^n$. Programs (\ref{eq:card-eig}) and (\ref{eq-variat-matrix-prog}) are equivalent, indeed if $X$ is a solution to the above problem, then $X \succeq 0$ and $\Rank(X)=1$ mean that we have $X=xx^T$, while $\Tr(X)=1$ implies that $\|x\|_2=1$. Finally, if $X=xx^T$ then $\Card(X) \leq k^2$ is equivalent to $\Card(x) \leq k$. We have made some progress by turning the convex maximization objective $x^T \Sigma x$ and the nonconvex constraint $\|x\|_2=1$ into a linear constraint and linear objective. Problem (\ref{eq-variat-matrix-prog}) is, however, still nonconvex and we need to relax both the rank and cardinality constraints. 

Since for every $u\in\reals^n$, $\Card(u)=q$ implies $\|u\|_1 \leq \sqrt{q} \|u\|_2$, we can replace the nonconvex constraint $\Card(X) \leq k^2$, by a weaker but convex constraint: $\ones ^T |X| \ones \leq k$, where we exploit the property that $\|X\|_F= \sqrt{x^Tx}=1$ when $X=xx^T$ and $\Tr(X)=1$. If we drop the rank constraint, we can form a relaxation of  (\ref{eq-variat-matrix-prog}) and (\ref{eq:card-eig}) as
\BEQ \label{eq-variat-relax} \BA{lll}
\mbox{maximize} & \Tr(\Sigma X)\\
\mbox{subject to} & \Tr(X)=1\\
& \ones ^T |X| \ones \leq k\\
& X \succeq 0,
\EA \EEQ
which is a semidefinite program in the variable $X\in\symm^n$, where $k$ is an integer parameter controlling the sparsity of the solution. The optimal value of this program will be an upper bound on the optimal value of the variational problem in (\ref{eq:card-eig}). Here, the relaxation of $\Card(X)$ in $\ones ^T |X| \ones$ corresponds to a classic technique which replaces the (non-convex) cardinality or $l_0$ norm of a vector $x$ with its largest convex lower bound on the unit box: $|x|$, the $l_1$ norm of $x$ (see \cite{Boyd00} or \cite{Dono05} for other applications).

Problem (\ref{eq-variat-relax}) can be interpreted as a robust formulation of the maximum eigenvalue problem, with additive, componentwise uncertainty in the input matrix $\Sigma$.  We again assume $A$ to be symmetric and positive semidefinite. If we consider a variation in which we penalize by the $\ell_1$ norm of the matrix $X$ instead of imposing a hard bound, to get
\BEQ
\label{eq:penalized-relax}
\BA{ll}
\mbox{maximize} & \Tr(\Sigma X) - \rho \ones^T |X| \ones\\
\mbox{subject to} & \Tr(X)=1\\
& X \succeq 0,
\EA \EEQ
which is a semidefinite program in the variable $X\in \symm^n$, where $\rho>0$ controls the magnitude of the penalty. We can rewrite this problem as
\BEQ
\label{eq:saddle}
\max_{X \succeq 0, \Tr(X)=1} ~ \min_{|U_{ij}|\leq \rho} \Tr(X(\Sigma +U)) 
\EEQ
in the variables $X\in\symm^n$ and $U\in\symm^n$. This yields the following dual to (\ref{eq:penalized-relax})
\BEQ
\label{eq:dual-robust}
\BA{ll}
\mbox{minimize} & \lambda^{\mathrm{max}}(\Sigma+U)\\
\mbox{subject to} & |U_{ij}|\leq \rho,\quad i,j=1,\ldots,n,
\EA \EEQ
which is a maximum eigenvalue problem with variable $U\in\symm^n$. This gives a natural robustness interpretation to the relaxation in (\ref{eq:penalized-relax}): it corresponds to a worst-case maximum eigenvalue computation, with componentwise bounded noise of intensity $\rho$ imposed on the matrix coefficients.

Finally, the KKT conditions (see \cite[\S5.9.2]{Boyd03}) for problems (\ref{eq:penalized-relax}) and (\ref{eq:dual-robust}) are given by
\BEQ
\label{eq:KKT}
\left\{
\BA{l}
(\Sigma+U)X=\lambda^{\mathrm{max}}(\Sigma+U)X\\
U\circ X=-\rho |X|\\
\Tr(X)=1,~X \succeq 0\\
|U_{ij}|\leq \rho,\quad i,j=1,\ldots,n.\\
\EA 
\right.
\EEQ
If the eigenvalue $\lambda^{\mathrm{max}}(\Sigma+U)$ is simple (when, for example, $\lambda^{\mathrm{max}}(A)$ is simple and $\rho$ is sufficiently small), the first condition means that $\Rank(X)=1$ and the semidefinite relaxation is \emph{tight}, with in particular $\Card(X)=\Card(x)^2$ if $x$ is the dominant eigenvector of $X$. When the optimal solution $X$ is not of rank one because of degeneracy (i.e. when $\lambda^{\mathrm{max}}(\Sigma+U)$ has multiplicity strictly larger than one), we can truncate $X$ as in \cite{Aliz95,Lema99}, retaining only the dominant eigenvector $x$ as an approximate solution to the original problem. In that degenerate scenario however, the dominant eigenvector of $X$ is not guaranteed to be as sparse as the matrix itself.

\subsection{A semidefinite relaxation with $\ell_0$ penalization} \label{ss:l0}
We summarize here the results in \cite{dAsp08b}. We begin by reformulating (\ref{eq:pca-card}) as a relatively simple convex maximization problem. Suppose that $\rho \geq \Sigma_{11}$. Since $z^T\Sigma z\leq \Sigma_{11}(\sum_{i=1}^n |z_i| )^2$ and $(\sum_{i=1}^n |z_i| )^2\leq \|z\|^2\Card(z)$ for all $z\in\reals^n$, we have
\[\BA{ll}
\phi(\rho) &= \max_{ \| z \|  \leq  1 } z^T \Sigma z - \rho \Card(z) \\
&\leq (\Sigma_{11}-\rho) \Card(z)\\  
&\leq 0,
\EA\]
hence the optimal solution to (\ref{eq:pca-card}) when $\rho \geq \Sigma_{11}$ is $z=0$. From now on, we assume $\rho \leq \Sigma_{11}$ in which case the inequality $\|z\|\leq1$ is tight. We can represent the sparsity pattern of a vector $z$ by a vector $u \in \{0,1\}^n$ and rewrite (\ref{eq:pca-card}) in the equivalent form
\[\BA{ll}
\phi(\rho) & = \dsp \max_{u \in \{0,1\}^n } \lambdamax(\diag(u)\Sigma\diag(u)) - \rho \ones^Tu\\
& = \dsp \max_{u \in \{0,1\}^n } \lambdamax(\diag(u)A^TA\diag(u)) - \rho \ones^Tu\\
& = \dsp \max_{u \in \{0,1\}^n } \lambdamax(A\diag(u)A^T) - \rho \ones^Tu,
\EA\]
using the fact that $\diag(u)^2=\diag(u)$ for all variables $u \in \{0,1\}^n $
and that for any matrix $B$, $\lambdamax(B^T B) = \lambdamax(BB^T)$. We then have
\[\BA{ll}
\phi(\rho) & = \dsp \max_{u \in \{0,1\}^n } \lambdamax(A\diag(u)A^T) - \rho \ones^Tu\\
& = \dsp \max_{ \| x \| = 1}~\max_{u \in \{0,1\}^n} x^TA\diag(u)A^Tx- \rho \ones^Tu\\
& = \dsp \max_{ \| x \| = 1}~\max_{u \in \{0,1\}^n}  \sum_{i=1}^n
u_i( (a_i^T x)^2 - \rho ).
\EA\]
Hence we finally get, after maximizing in $u$ (and using $\max_{v \in \{0,1\}} \beta v = \beta_+$)
\BEQ \label{eq:pca-ncvx}
\phi(\rho) = \max_{ \| x \| = 1} \sum_{i=1}^n((a_i^Tx)^2-\rho)_+,
\EEQ
which is a nonconvex problem in the variable $x\in\reals^n$. We then select variables $i$ such that $(a_i^Tx)^2-\rho>0$. Note that if $\Sigma_{ii} = a_i^Ta_i < \rho$, we must have $(a_i^Tx)^2 \leq \|a_i\|^2 \|x\|^2 < \rho$ hence variable $i$ will never be part of the optimal subset and we can remove it.

Because the variable $x$ appears solely through $X=xx^T$, we can reformulate the problem in terms of $X$ only, using the fact that when $\|x\|=1$, $X=xx^T$ is equivalent to $\Tr(X)=1$, $X\succeq 0$ and $\Rank(X)=1$. We thus  rewrite (\ref{eq:pca-ncvx}) as
\[\BA{lll}
\phi(\rho) =&\mbox{max.}&  \sum_{i=1}^n(a_i^TXa_i-\rho)_+ \\
&\mbox{s.t.} & \Tr(X)=1,~\Rank(X)=1\\
& & X\succeq 0.\\
\EA\]
Note that because we are maximizing a convex function $\Delta_n=\{X\in\symm_n:~\Tr(X)=1,~ X\succeq 0\}$ which is convex, the solution must be an extreme point of $\Delta_n$ (i.e. a rank one matrix), hence we can drop the rank constraint here. Unfortunately, $X \mapsto (a_i^TXa_i-\rho)_+$, the function we are \emph{maximizing}, is convex in $X$ and not concave, which means that the above problem is still hard. However, we show below that on rank one elements of $\Delta_n$, it is also equal to a concave function of $X$, and we use this to produce a semidefinite relaxation of problem (\ref{eq:pca-card}).

\begin{proposition} \label{prop:relax} Let $A\in\reals^{n \times n}$, $\rho \geq 0$ and denote by $a_1,\dots,a_n \in \reals^n$ the columns of $A$, an upper bound on
\BEQ\label{eq:pca-plus}\BA{lll}
\phi(\rho) =&\mbox{max.}&  \sum_{i=1}^n(a_i^TXa_i-\rho)_+ \\
&\mbox{s.t.} & \Tr(X)=1,~ X\succeq 0, ~\Rank(X)=1\EA\EEQ
can be computed by solving
\BEQ  \label{eq:pca-relax}
\BA{lll}
\psi(\rho) =&\mbox{max.}&  \sum_{i=1}^n \Tr(X^{1/2}B_i X^{1/2})_+ \\
&\mbox{s.t.} & \Tr(X)=1,~X\succeq 0.
\EA\EEQ
in the variables $X\in\symm_n$, where $B_i=a_ia_i^T - \rho\idm$, or also
\BEQ \label{eq:primal-p}
\BA{lll}
\psi(\rho) = &\mbox{max.}&  \sum_{i=1}^n\Tr(P_iB_i) \\
&\mbox{s.t.} & \Tr(X)=1,~X\succeq 0,~X\succeq P_i \succeq 0,
\EA\EEQ 
which is a semidefinite program in the variables $X\in\symm_n,~P_i\in\symm_n$.
\end{proposition}
\begin{proof} We let $X^{1/2}$ be the positive square root (i.e. with nonnegative eigenvalues) of a symmetric positive semi-definite matrix $X$. In particular, if $X=xx^T$ with $\|x\|=1$, then $X^{1/2}= X = xx^T$, and for all  $\beta\in \reals$, $\beta xx^T$ has one eigenvalue equal to $\beta$ and $n-1$ equal to 0, which implies $\Tr(\beta xx^T)_+=\beta_+$. We thus get
\BEAS
 (a_i^TXa_i-\rho)_+ & = & \Tr ( (a_i^T xx^T a_i - \rho) xx^T )_+
 \\
& = &  \Tr(x ( x^T a_ia_i^T x -\rho ) x^T )_+ \\
& = &  \Tr( X^{1/2} a_i a_i^T X^{1/2} - \rho X  )_+
= \Tr( X^{1/2} (  a_i a_i^T -\rho \idm)  X^{1/2}    )_+.
\EEAS
For any symmetric matrix $B$, the function $X \mapsto \Tr( X^{1/2} B X^{1/2} )_+$ is concave on the set of symmetric positive semidefinite matrices, because we can write it as
\[\BA{ll}
\Tr(X^{1/2}BX^{1/2})_+&=\dsp\max_{\{0\preceq P \preceq X\}} \Tr(PB)\\
&=\dsp\min_{\{Y\succeq B,~Y\succeq 0\}} \Tr(YX),
\EA\]
where this last expression is a concave function of $X$ as a pointwise minimum of affine functions. We can now relax the original problem into a convex optimization problem by simply dropping the rank constraint, to get
$$ \BA{lll}
\psi(\rho) \equiv &\mbox{max.}&  \sum_{i=1}^n\Tr(X^{1/2}a_ia_i^TX^{1/2}-\rho X)_+ \\
&\mbox{s.t.} & \Tr(X)=1,~X\succeq 0,
\EA
$$
which is a convex program in $X\in\symm_n$. Note that because $B_i$   has at most one nonnegative eigenvalue, we can replace $\Tr(X^{1/2}a_ia_i^TX^{1/2}-\rho X)_+$ by $\lambdamax(X^{1/2}a_ia_i^TX^{1/2}-\rho X)_+$ in the above program. Using the representation of $\Tr(X^{1/2}BX^{1/2})_+$ detailed above, problem (\ref{eq:pca-relax}) can be written as a semidefinite program
\[
\BA{lll}
\psi(\rho) = &\mbox{max.}&  \sum_{i=1}^n\Tr(P_iB_i) \\
&\mbox{s.t.} & \Tr(X)=1,~X\succeq 0,~X\succeq P_i \succeq 0,
\EA\]
in the variables $X\in\symm_n,~P_i\in\symm_n$, which is the desired result. 
\end{proof}

Note that we always have $\psi(\rho)\geq\phi(\rho)$ and when the solution to the above semidefinite program has rank one, $\psi(\rho)=\phi(\rho)$ and the semidefinite relaxation (\ref{eq:primal-p}) is \emph{tight}. This simple fact allows us to derive sufficient global optimality conditions for the original sparse PCA problem. We recall in particular the following result from \cite{dAsp08b} which provides sufficient conditions for a particular nonzero coefficient pattern~$I$ to be globally optimal. The optimal solution $x$ to (\ref{eq:pca-card}) is then found by solving an eigenvalue problem on the principal submatrix of $\Sigma$ with support~$I$. 

\begin{proposition} \label{th:opt}
Let $A\in\reals^{n \times n}$, $\rho \geq 0$, $\Sigma=A^TA$ with $a_1,\dots,a_n \in \reals^n$ the columns of $A$. Given a sparsity pattern $I$, setting $x$ to be the largest eigenvector of $\sum_{i \in I} a_i a_i^T$, if there is a $\rho^*\geq 0$ such that the following conditions hold
\[
\max_{i \in I^c} (a_i^T x)^2 < \rho^* < \min_{i \in I} (a_i^T x)^2 \quad \mbox{and}\quad
\lambdamax\left(\sum_{i=1}^n Y_i\right) \leq \sum_{i \in I}((a_i^T x)^2 - \rho^*),
\] 
with the dual variables $Y_i$ defined as 
\[
Y_i=\max\left\{0,\rho\frac{(a_i^T a_i-\rho)}{(\rho-(a_i^T x)^2)}\right\}\frac{(\idm-xx^T)a_ia_i^T(\idm-xx^T)}{\|(\idm-xx^T)a_i\|^2} ,\quad{when }~ i \in I^c,
\]
and
\[
Y_i = \frac{ B_i x x^T B_i }{x^T B_i x},\quad{when }~ i \in I,
\]
then the sparsity pattern $I$ is globally optimal for the sparse PCA problem~(\ref{eq:pca-card}) with $\rho=\rho^*$ and we can form an optimal solution $z$ by solving the maximum eigenvalue problem
\[
z=\argmax_{ \{z_{I^c}=0,~\|z\|=1 \}} z^T \Sigma z.\\
\]
\end{proposition}
This result also provides tractable {\em lower bounds} on the optimal value of~(\ref{eq:pca-card}) whenever the solution is not optimal.

\section{Algorithms} \label{s:algos}
In this section, we describe algorithms for solving the semidefinite relaxations detailed above. We also describe greedy methods to improve the quality of these solutions.

\subsection{First-order methods} \label{s:smooth}
Again, given a covariance matrix $\Sigma\in\symm_n$, the DSPCA code solves a penalized formulation of problem~(\ref{eq-variat-relax}), written as
\BEQ\label{eq:orig-pb} \BA{ll}
\mbox{maximize} & \Tr(\Sigma X) - \rho \ones^T |X| \ones\\
\mbox{subject to} & \Tr(X)=1\\
& X \succeq 0,
\EA
\EEQ
in the variable $X\in\symm_n$.  The dual of this program can be written as
\BEQ\label{eq:orig-dual}
\BA{ll}
\mbox{minimize} & f(U)=\lambda^\mathrm{max}(\Sigma +U)\\
\mbox{subject to} & |U_{ij}|\leq \rho.\\
\EA
\EEQ
in the variable $U\in\symm^n$. The algorithm in \cite{dasp04a,Nest04a} regularizes the objective $f(U)$ in (\ref{eq:orig-dual}), replacing it by the smooth (i.e. with Lipschitz continuous gradient) uniform approximation
\[
f_\mu(U) = \mu \log \left(\Tr \exp((\Sigma+U)/\mu)\right) - \mu \log n.
\]
Following \cite{nest83}, solving the smooth problem
\[
\min_{U\in{\cal Q}} f_\mu(U)
\]
where ${\cal Q}=\{U\in{\footnotesize \symm^n},\,|U_{ij}|\leq \rho\}$, with $\mu =\epsilon/2\log(n)$ then produces an $\epsilon$-approximate solution to (\ref{eq:orig-pb}). The key difference between the minimization scheme developed in \cite{nest83} and classical gradient minimization methods is that it is not a descent method but achieves a complexity of $O(L/N^2)$ instead of $O(1/N)$ for gradient descent, where $N$ is the number of iterations and $L$ the Lipschitz constant of the gradient. Furthermore, this convergence rate is provably optimal for this particular class of convex minimization problems (see \cite[Th.~2.1.13]{Nest03a}). Thus, by sacrificing the (local) properties of descent directions, we improve the (global) complexity estimate by an order of magnitude. For our problem here, once the regularization parameter $\mu$ is set, the algorithm is detailed as Algorithm~\ref{alg:first-order}.

\begin{algorithm}[h]
\caption{First-Order Algorithm.}
\label{alg:first-order} 
\begin{algorithmic} [1]
\REQUIRE The covariance $\Sigma \in \reals^{n \times n}$, and a parameter $\rho>0$ controlling sparsity.
\FOR{$i=1$ to $N$} 
\STATE Compute $f_\mu(U_i)$ and $\nabla f_\mu(U_i)$
\STATE Find $Y_i = \arg\min_{Y \in {\cal Q}} \: \langle \nabla f_\mu(U_i) , Y \rangle + \frac{1}{2} L \|U_i-Y\|_F^2$
\STATE Find $W_i = \arg\min_{W\in {\cal Q}} \left\{ \frac{L \|W\|_F^2}{2}  + \sum_{j=0}^N \frac{j+1}{2}(f_\mu(U_j)+\langle \nabla f_\mu(U_j),W-U_j \rangle ) \right\}$
\STATE Set $U_{i+1} = \frac{2}{i+3} W_i + \frac{i+1}{i+3} Y_i$
\ENDFOR 
\ENSURE A matrix $U\in\symm_n$.
\end{algorithmic} 
\end{algorithm} 

The algorithm has four main steps. Step one computes the (smooth) function value and gradient. The second step computes the \emph{gradient mapping}, which matches the gradient step for unconstrained problems (see \cite[p.86]{Nest03a}). Step three and four update an \emph{estimate sequence} see (\cite[p.72]{Nest03a}) of $f_\mu$ whose minimum can be computed explicitly and gives an increasingly tight upper bound on the minimum of $f_\mu$. We now present these steps in detail for our problem (we write $U$ for $U_i$ and $X$ for $X_i$).

\paragraph{Step 1.} The most expensive step in the algorithm is the first, the computation of $f_\mu$ and its gradient. The function value can be reliably computed as
\[
f_\mu(U) = d_{\rm max} + \mu \log \left( \sum_{i=1}^n \exp(\frac{d_i-d_{\rm max}}{\mu}) \right) - \mu \log n .
\]
where $d_i$ are the eigenvalues of $\Sigma+U$. The gradient $\nabla f_\mu(U)$ can be computed explicitly as
\[
\nabla f_\mu(U):=  \exp\left((\Sigma+U)/\mu \right)/\Tr\left(\exp\left((\Sigma +U)/\mu \right)\right).
\]
which means computing the same matrix exponential. 

\paragraph{Step 2.}
This step involves a problem of the form
\[
\arg\min_{Y \in {\cal Q}} \: \langle \nabla f_\mu(U) , Y \rangle + \frac{1}{2} L \|U-Y\|_F^2 ,
\]
where $U$ is given.  The above problem can be reduced to a Euclidean projection
\begin{equation}\label{eq:euclo-prox}
\arg\min_{\|Y\|_\infty \leq 1} \: \|Y - V\|_F,
\end{equation}
where $V = U - L^{-1}\nabla f_\mu(U)$ is given. The solution is given by
\[
Y_{ij} = \mbox{\bf sgn}(V_{ij}) \min(|V_{ij}|,1), \quad i,j=1,\ldots,n.
\]

\paragraph{Step 3.} The third step involves solving a Euclidean projection problem similar to (\ref{eq:euclo-prox}), with the solution $V$ defined by
\[
V = -\frac{1}{L} \sum_{i=0}^k \frac{i+1}{2} \nabla f_\mu(U_i).
\]

\paragraph{Stopping criterion} We can stop the algorithm when the duality gap is smaller than $\epsilon$:
\[
\mathrm{gap}_k=\lambdamax(\Sigma+U_k) - \Tr \Sigma X_i + \ones^T |X_i|\ones \leq \epsilon,
\]
where $X_k = \nabla f_\mu(U)$ is our current estimate of the dual variable. The above gap is necessarily non-negative, since both $X_i$ and $U_i$ are feasible for the primal and dual problem, respectively. This is checked periodically, for example every $100$ iterations. 

\paragraph{Complexity} Overall the algorithm requires 
\BEQ\label{eq:complexity}
O\left(\rho \frac{n\sqrt{\log n}}{\epsilon}\right)
\EEQ
iterations \cite{dasp04a,Nest04a}. The main step at each iteration is computing the matrix exponential $\exp((\Sigma+U)/\mu)$ (see \cite{Mole03} for a comprehensive survey) at a cost of $O(n^3)$ flops.

\subsection{Greedy methods} \label{s:greedy}
\label{sec:greedy} We can also find good solution to problem (\ref{eq:pca-card}), or improve existing solutions, using greedy methods. We first present very simple preprocessing solutions with complexity $O(n\log n)$ and $O(n^2)$. We then recall a simple greedy algorithm with complexity $O(n^4)$. Finally, our first contribution in this section is to derive an approximate greedy algorithm that computes a full set of (approximate) solutions for problem (\ref{eq:pca-card}), with complexity $O(n^3)$.

\subsubsection{Sorting and thresholding}
The simplest ranking algorithm is to sort the diagonal of the matrix $\Sigma$ and rank the variables by variance. This works intuitively because the diagonal is a rough proxy for the eigenvalues: the Schur-Horn theorem states that the diagonal of a matrix majorizes its  eigenvalues \cite{Horn85}; sorting costs $O(n\log n)$. Another quick solution is to compute the leading eigenvector of $\Sigma$ and form a sparse vector by thresholding to zero the coefficients whose magnitude is smaller than a certain level. This can be done with cost $O(n^2)$.

\subsubsection{Full greedy solution}
\label{sec:greedy-classic} Following \cite{Mogh06b}, starting from an initial solution of cardinality one at $\rho=\Sigma_{11}$, we can update an increasing sequence of index sets $I_k\subseteq[1,n]$, scanning all the remaining variables to find the index with maximum variance contribution. 

\begin{algorithm} 
\caption{Greedy Search Algorithm.} 
\label{alg:greedy-search} 
\begin{algorithmic} [1]
\REQUIRE $\Sigma \in \reals^{n \times n}$
\STATE Preprocessing: sort variables by decreasing diagonal elements and permute elements of $\Sigma$ accordingly. 
\STATE Compute the Cholesky decomposition $\Sigma = A^T A $.
\STATE Initialization: $I_1 = \{1\}$,  $x_1=a_1/\|a_1\|$.
\FOR{$i=1$ to $k^\mathrm{target}$} 
\STATE Compute $i_k = \argmax_{i \notin I_k}  \lambdamax\left(
\sum_{j\in I_k \cup \{ i \} } a_j a_j^T \right)$.
\STATE Set $I_{k+1}=I_k\cup\{i_k\}$ and compute $x_{k+1}$ as the leading eigenvector of $\sum_{ j \in I_{k+1}}a_ja_j^T$.
\ENDFOR 
\ENSURE Sparsity patterns $I_k$.
\end{algorithmic} 
\end{algorithm} 

At every step, $I_k$ represents the set of nonzero elements (or sparsity pattern) of the current point and we can define $z_k$ as the solution to problem (\ref{eq:pca-card}) given $I_k$, which is:
\[
z_k=\argmax_{ \{z_{I_k^c}=0,~\|z\|=1 \}} z^T \Sigma z - \rho k,\\
\]
which means that $z_k$ is formed by padding zeros to the leading eigenvector of the submatrix $\Sigma_{I_k,I_k}$. Note that the entire algorithm can be written in terms of a factorization $\Sigma = A^T A$ of the matrix $\Sigma$, which means significant computational savings when $\Sigma$ is given as a Gram matrix. The matrices $\Sigma_{I_k,I_k}$ and $\sum_{i \in I_k} a_i a_i^T$ have the same eigenvalues and if $z$ is an eigenvector of $\Sigma_{I_k,I_k}$, then $ A_{I_k} z / \|A_{I_k} z \|$ is an eigenvector of $A_{I_k} A_{I_k}^T$.

\subsubsection{Approximate greedy solution}
\label{sec:greedy-path} Computing $n-k$ eigenvalues at each iteration is costly and we can use the fact that $uu^T$ is a subgradient of $\lambdamax$ at $X$ if $u$ is a leading eigenvector of $X$~\cite{Boyd03}, to get: 
\BEQ
\label{eq:lambdamax-lower-bound}
\lambdamax\left(
\sum_{j\in I_k \cup \{ i \} } a_j a_j^T \right)
\geq 
\lambdamax\left(
\sum_{j\in I_k } a_j a_j^T \right)
+ (x_k^T a_i)^2,
\EEQ
which means that the variance is increasing by at least $(x_k^T a_i)^2$ when variable $i$ is added to $I_k$. This provides  a lower bound on the objective which does not require finding $n-k$ eigenvalues at each iteration. We then derive the following algorithm.

\begin{algorithm} 
\caption{Approximate Greedy Search Algorithm.} 
\label{alg:approx-greedy-search} 
\begin{algorithmic} [1]
\REQUIRE $\Sigma \in \reals^{n \times n}$
\STATE Preprocessing: sort variables by decreasing diagonal elements and permute elements of $\Sigma$ accordingly. 
\STATE Compute the Cholesky decomposition $\Sigma = A^T A $.
\STATE Initialization: $I_1 = \{1\}$,  $x_1=a_1/\|a_1\|$.
\FOR{$i=1$ to $k^\mathrm{target}$} 
\STATE Compute $i_k = \argmax_{i \notin I_k}  (x_k^T a_i)^2 $.
\STATE Set $I_{k+1}=I_k\cup\{i_k\}$ and compute $x_{k+1}$ as the leading eigenvector of $\sum_{ j \in I_{k+1}}a_ja_j^T$.
\ENDFOR 
\ENSURE Sparsity patterns $I_k$.
\end{algorithmic} 
\end{algorithm} 

Again, at every step, $I_k$ represents the set of nonzero elements (or sparsity pattern) of the current point and we can define $z_k$ as the solution to problem (\ref{eq:pca-card}) given $I_k$, which is:
\[
z_k=\argmax_{ \{z_{I_k^c}=0,~\|z\|=1 \}} z^T \Sigma z - \rho k,\\
\]
which means that $z_k$ is formed by padding zeros to the leading eigenvector of the submatrix $\Sigma_{I_k,I_k}$.  Better points can be found by testing the variables corresponding to the $p$ largest values of $(x_k^T a_i)^2$ instead of picking only the best one.

\subsubsection{Computational complexity}
The complexity of computing a greedy regularization path using the classic greedy algorithm in \S\ref{sec:greedy-classic} is $O(n^4)$: at each step $k$, it computes $(n-k)$ maximum eigenvalue of matrices with size $k$. The approximate algorithm in \S\ref{sec:greedy-path} computes a full path in $O(n^3)$: the first Cholesky decomposition is $O(n^3)$, while the complexity of the $k$-th iteration is $O(k^2)$ for the maximum eigenvalue problem and $O(n^2)$ for computing all products $(x^T a_j)$. Also, when the matrix $\Sigma$ is directly given as a Gram matrix $A^T A$ with $A\in\reals^{q \times n}$ with $q<n$, it is advantageous to use $A$ directly as the square root of $\Sigma$ and the total complexity of getting the path up to cardinality $p$ is then reduced to $O(p^3 + p^2 n)$ (which is $O(p^3)$ for the eigenvalue problems and $O(p^2 n )$ for computing the vector products).

\section{Applications} \label{s:apps}
In this section, we illustrate the sparse PCA approach in applications: in news (text), finance and voting data from the US Senate.

\subsection{News data}
Our first dataset is a small version of the ``20-newsgroups'' data\footnote{available from {\tt http://cs.nyu.edu/\~{}roweis/data.html}.}. The data records binary occurrences of $100$ specific words across $16242$ postings, where the postings have been tagged by the highest level domain in Usenet. 

\begin{figure}[htb]
\begin{center}
\begin{tabular}{cc}
\psfrag{PCA}{}
\psfrag{number of principal components}[t][b]{Number of Principal Comp.}
\psfrag{cumulative percentage of variance explained}[b][t]{\% Variance Explained}
\includegraphics[width = 0.48\textwidth,height = 0.38\textwidth]{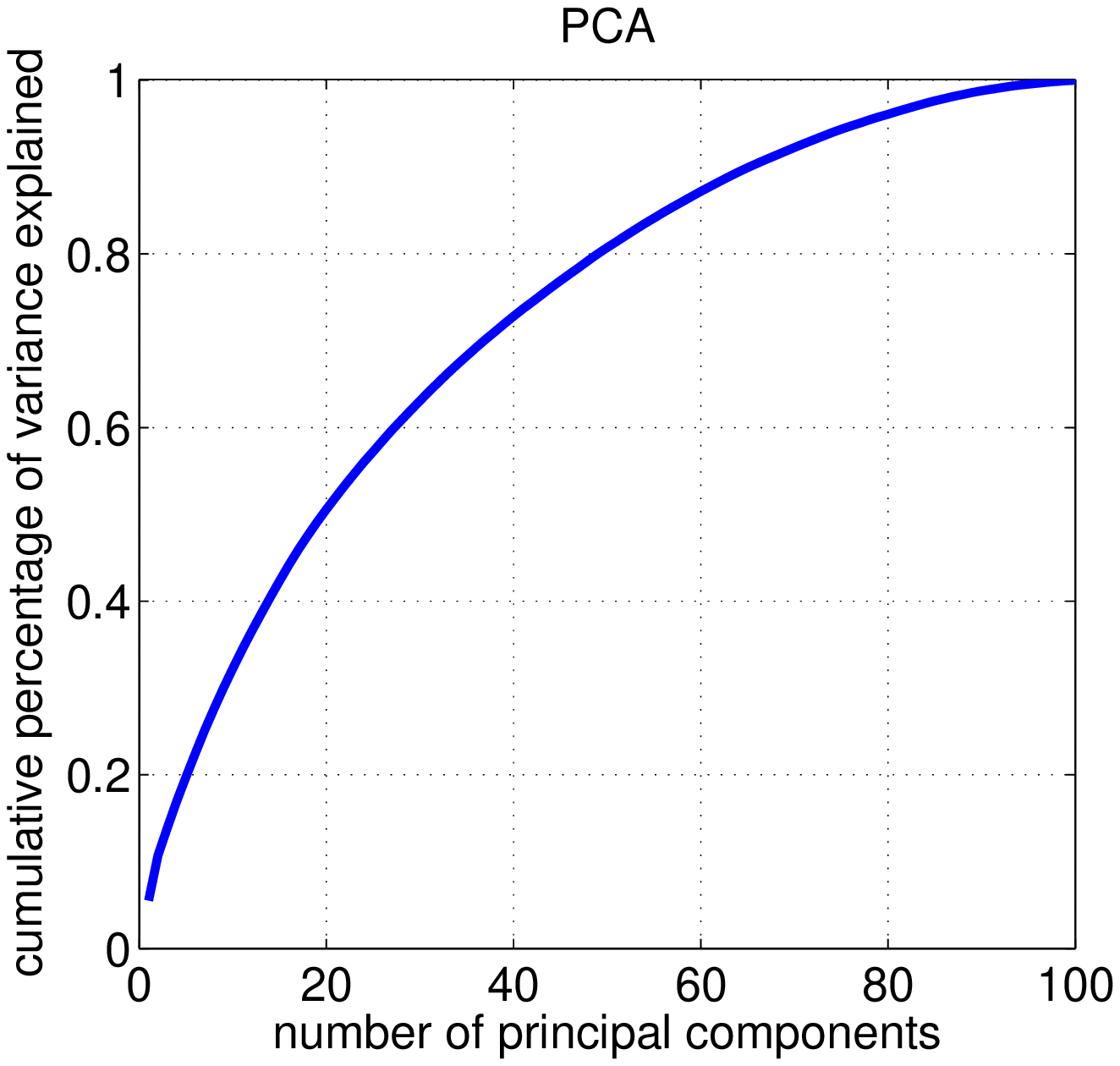} &
\psfrag{PCA}{}
\psfrag{2nd PC (cardinality = 100 )}[c][c]{\footnotesize{2nd P.C.}}
\psfrag{1st PC (cardinality = 100 )}[c][c]{\footnotesize{1st P.C.}}
\psfrag{3rd PC (cardinality = 100 )}[b][t]{\footnotesize{3rd P.C.}}
\includegraphics[width = 0.48\textwidth,height = 0.38\textwidth]{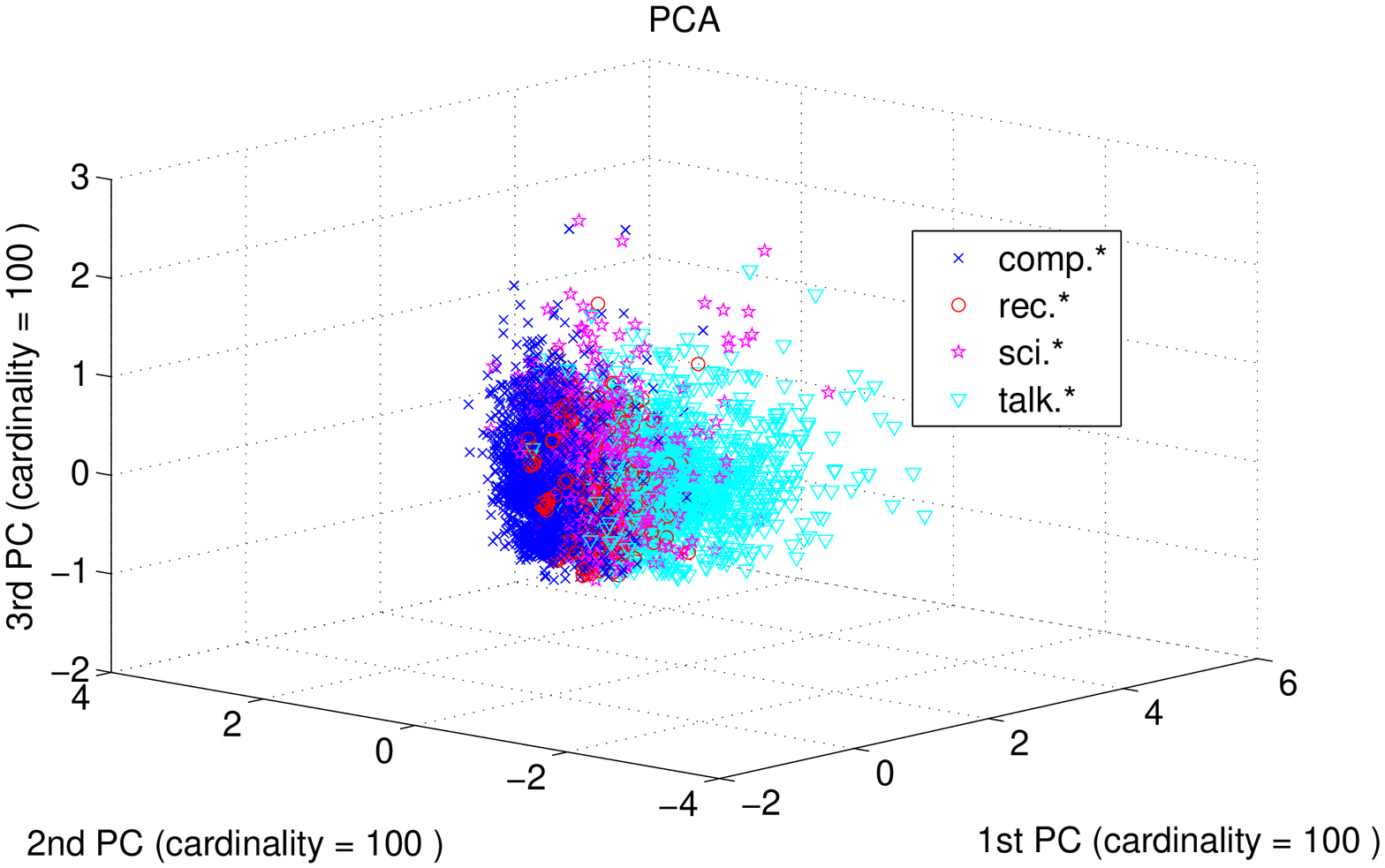}
\end{tabular}
\caption{PCA with 20-Newsgroups data. {\em Left:} Explained variance vs.\ number of PCs. {\em Right:} 3D visualization via PCA. \label{fig:ngVar}}
\end{center}
\end{figure}

Each posting is viewed as one point in a $100$-dimensional space. We begin with a standard PCA on the data. Fig.~\ref{fig:ngVar} (left) shows the cumulative percentage of variance explained as we increase the number of principal components. The slow increase means that the data does not lie within a subspace of significantly low dimension. We can anyway proceed to visualize the data: Fig.~\ref{fig:ngVar} (right) is the result obtained by projecting it on a subspace of dimension $3$, chosen by selecting the eigenvectors corresponding to the three largest eigenvalues of the covariance matrix. Since these $3$ vectors are dense, the axes in Fig.~\ref{fig:ngVar} (right) do not have a clear interpretation.

\begin{figure}[htb]
\begin{center}
\begin{tabular}{cc}
\psfrag{Cardinality}[t][b]{Cardinality}
\psfrag{Variance}[b][t]{\% Variance Explained}
\psfrag{DSPCA --- 1st Principal Component}[b][c]{1st P.C.}
\includegraphics[width = 0.48\textwidth, height = 0.38\textwidth]{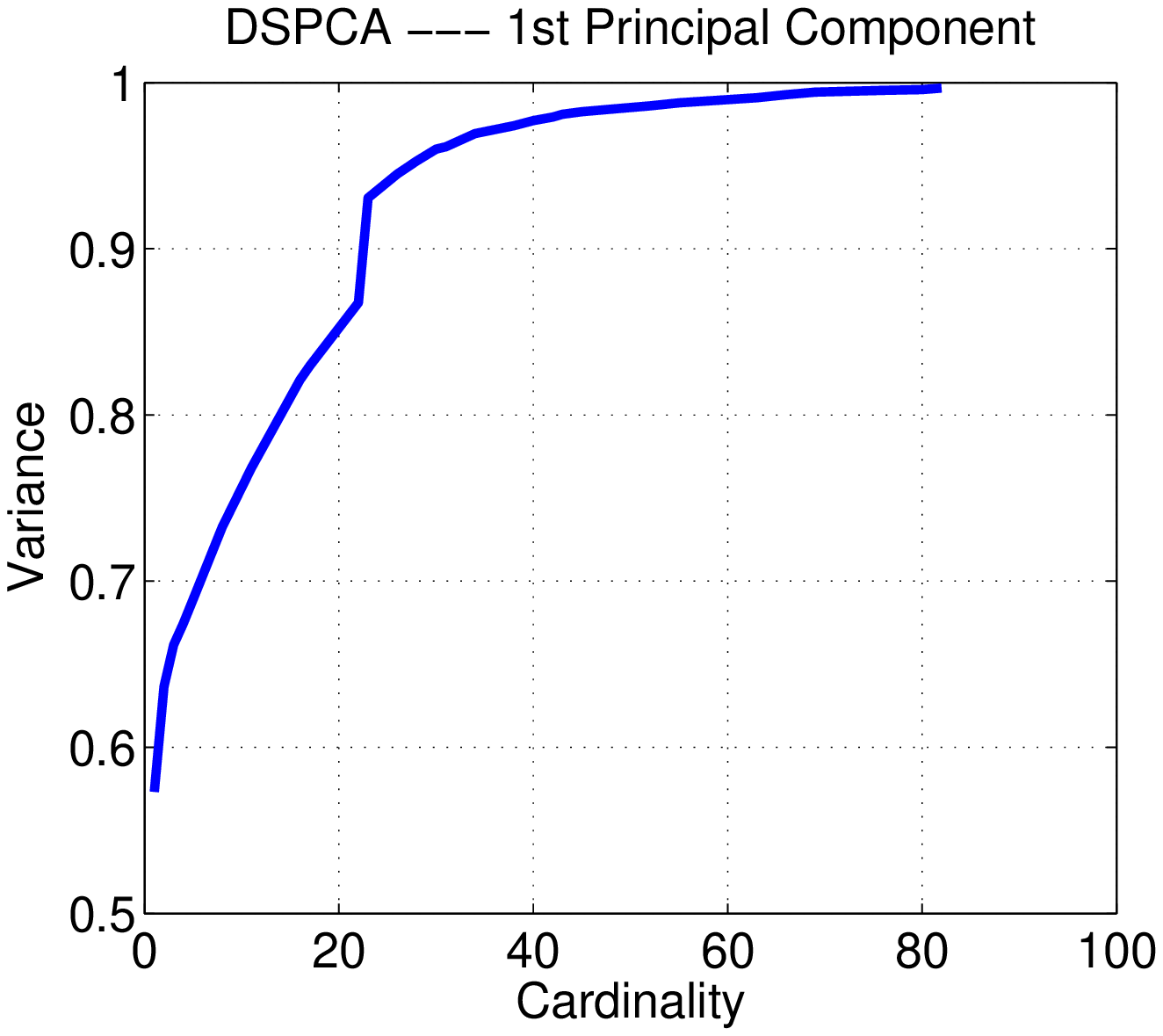}&
\psfrag{Cardinality}[t][b]{Cardinality}
\psfrag{Variance}[b][t]{\% Variance Explained}
\psfrag{DSPCA --- 2nd Principal Component}[b][c]{1st P.C.}
\includegraphics[width = 0.48\textwidth, height = 0.38\textwidth]{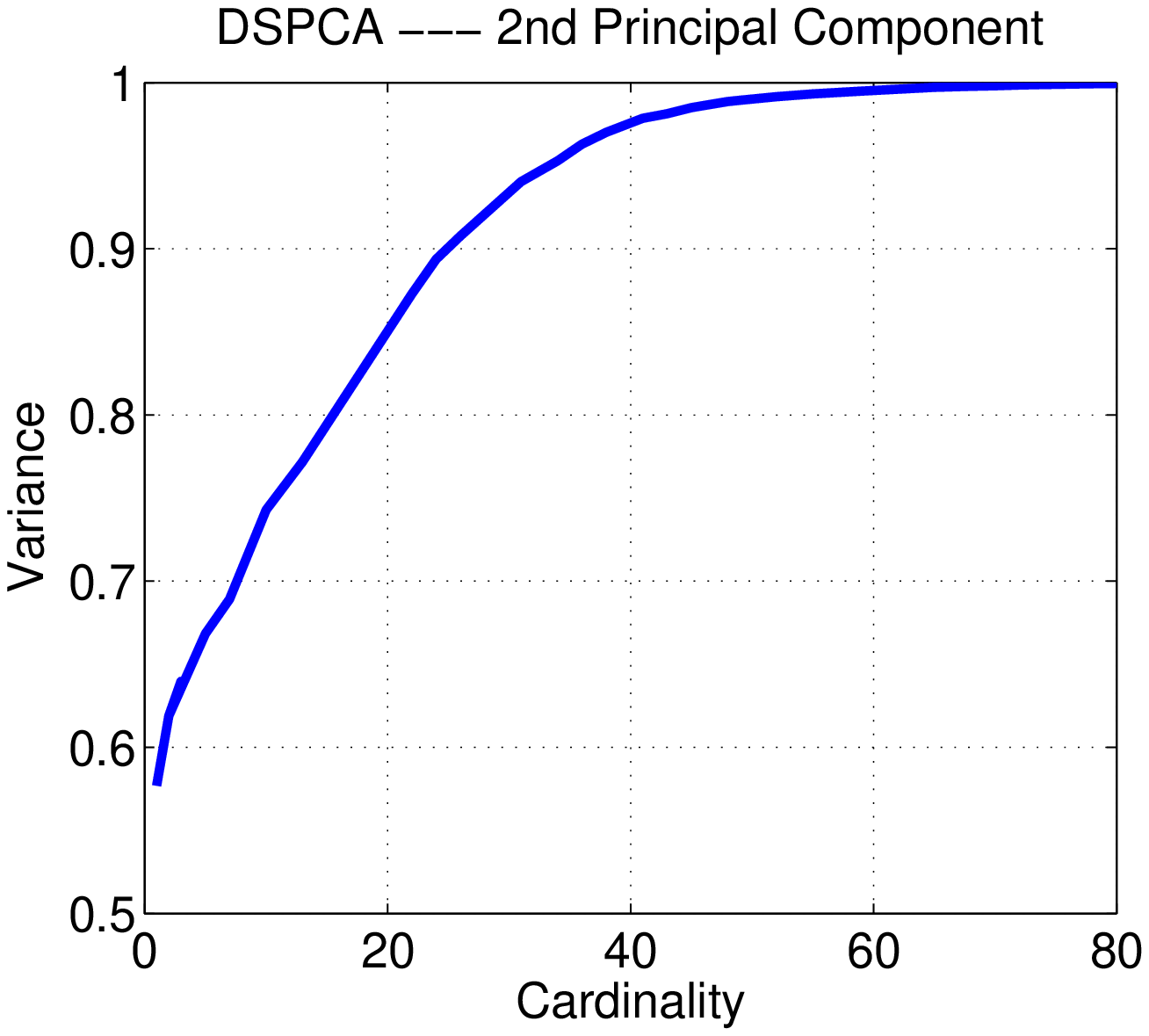}\\
\\
\\
\psfrag{Cardinality}[t][b]{Cardinality}
\psfrag{Variance}[b][t]{\% Variance Explained}
\psfrag{DSPCA --- 3rd Principal Component}[b][c]{1st P.C.}
\includegraphics[width = 0.48\textwidth, height = 0.38\textwidth]{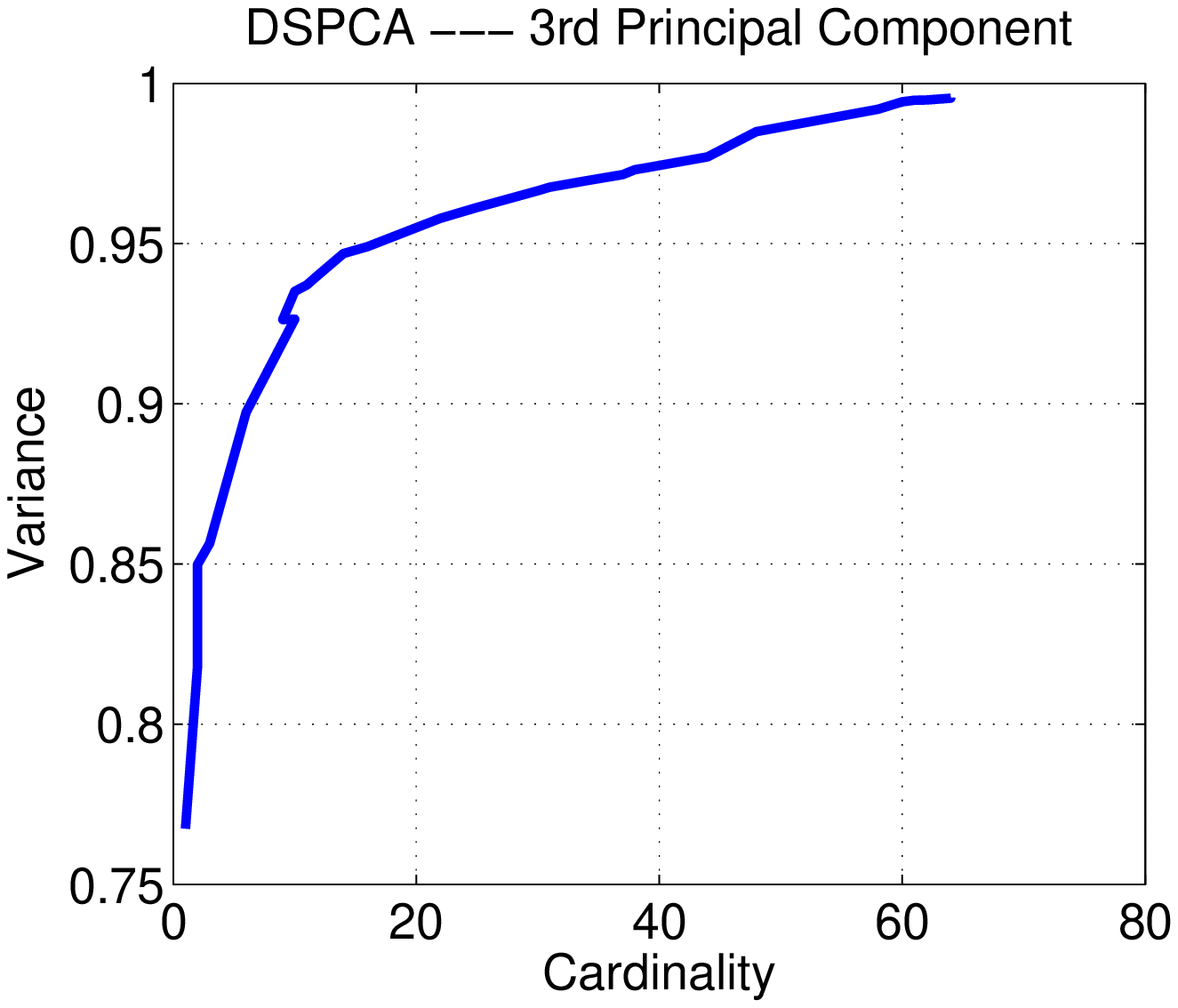}&
\psfrag{2nd PC (cardinality = 26 )}[c][c]{\footnotesize{2nd P.C.}}
\psfrag{1st PC (cardinality = 30 )}[c][c]{\footnotesize{1st P.C.}}
\psfrag{3rd PC (cardinality = 10 )}[b][t]{\footnotesize{3rd P.C.}}
\psfrag{Sparse PCA --- DSPCA}[b][c]{Sparse PCA}
\includegraphics[width = 0.48\textwidth, height = 0.38\textwidth]{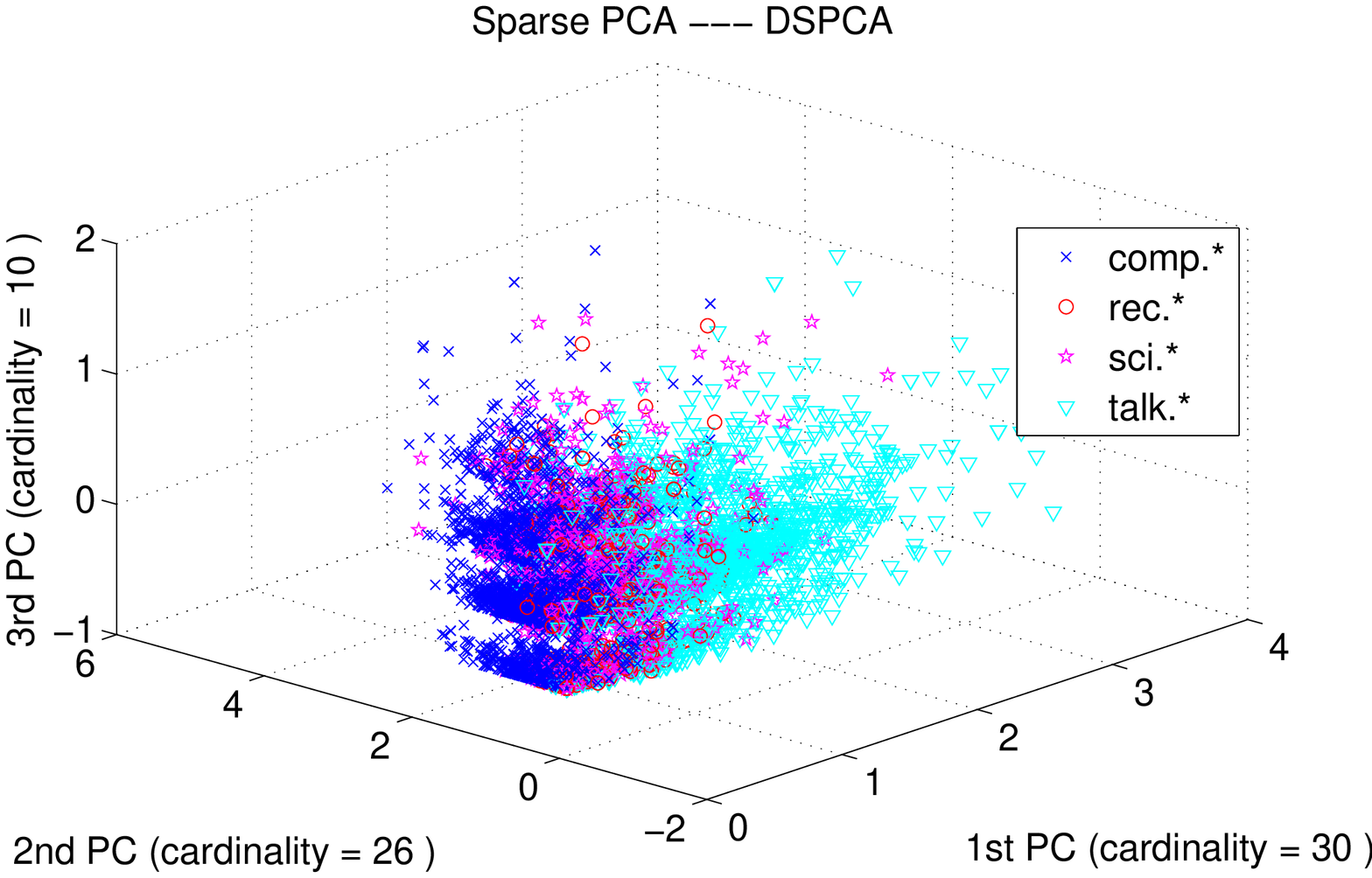}
\end{tabular}
\end{center}
\caption{Sparse PCA on the 20Newsgroups data set. First three principal components  and 3D visualization. The first three principal components have cardinalities 26, 30 and 10 respectively.\label{fig:ng1PC}}
\end{figure}

With sparse PCA, we hope to find a set of corresponding sparse principal components, which still help with visualization nearly as well as PCA does, and yet reveal some interesting structure. To achieve this, we have run the first-order algorithm of Section~\ref{s:smooth} (referred to as "DSPCA" hereafter) on the data with a range of values for the penalty parameter $\rho$. We obtained a plot of the variance explained by the first sparse principal component (PC), as a function of its cardinality (Fig.~\ref{fig:ng1PC}). We then selected a cardinality that can explain at least $90\%$ of the variance explained by the the first principal component obtained from PCA. Then we have deflated the covariance matrix by taking out the part due to the first sparse PC, and then repeated the above procedure to obtain the second sparse PC. In the same way, we have solved for the third sparse PC. Fig.~\ref{fig:ng1PC} also shows the projection of the data on the $3$-dimensional subspace that is spanned by the three sparse PCs obtained above. 

We first note that only a small number of words, out of the total of $100$ words that can appear in each sparse PC, can explain more than $90\%$ of variance explained by the corresponding PC. Specifically, we obtain $30$ words for the first PC, $26$ for the second, and $10$ for the third.
The lists of words associated with each sparse PCs is given in Table \ref{tab:ngloading}, and reveals some structure about each one of the sparse PCs. That is, the $30$ words associated with the first sparse PC are almost all about politics and religion, the $26$ words in the second sparse PC are all computer-related, and the majority of the $10$ words in the third sparse PC concerns science.  Hence, applying sparse PCA to this data set allows to discover structure that is otherwise hidden in the standard PCA, for example that the first principal component is mainly related to politics and religion.

{\small \begin{table}[htb]
\centering \small{\texttt{
\begin{tabular}{ccc}
 {\rm ~1st PC (30 words)~} & {\rm ~2nd PC (26 words)~} & {\rm ~3rd PC (10 words)~} \\
\hline
\texttt{fact} & \texttt{help} & \texttt{problem} \\
\texttt{question} & \texttt{problem} & \texttt{university} \\
    \texttt{world} & \texttt{system} & \texttt{email} \\
    \texttt{course} & \texttt{email}  & \texttt{state} \\
    case & windows & research \\
    problem & program& science \\
    god & computer & phone \\
    government & software & world \\
    human & university & fact \\
    state & version & question \\
    number & files &\\
    christian & drive &\\
    evidence & data &\\
    law & card &\\
    power & dos&\\
    religion & god&\\
    children & disk&\\
    jesus & pc&\\
    system & graphics&\\
    rights & ftp&\\
    war & memory&\\
    jews & christian&\\
    help & phone&\\
    bible & video&\\
    earth & fact&\\
    science & display&\\
    research & & \\
    israel & & \\
    president &&\\
    gun &&\\
 \hline
\end{tabular}}}
\caption{\label{tab:ngloading}
Words associated with the first three sparse PCs.
}
\end{table}}

\begin{figure}[htb]
\begin{center}
\begin{tabular}{cc}
\psfrag{Cardinality}[t][b]{Cardinality}
\psfrag{\% Variance Explained}[b][t]{\% Variance Explained}
\includegraphics[width = 0.48\textwidth, height = 0.38\textwidth]{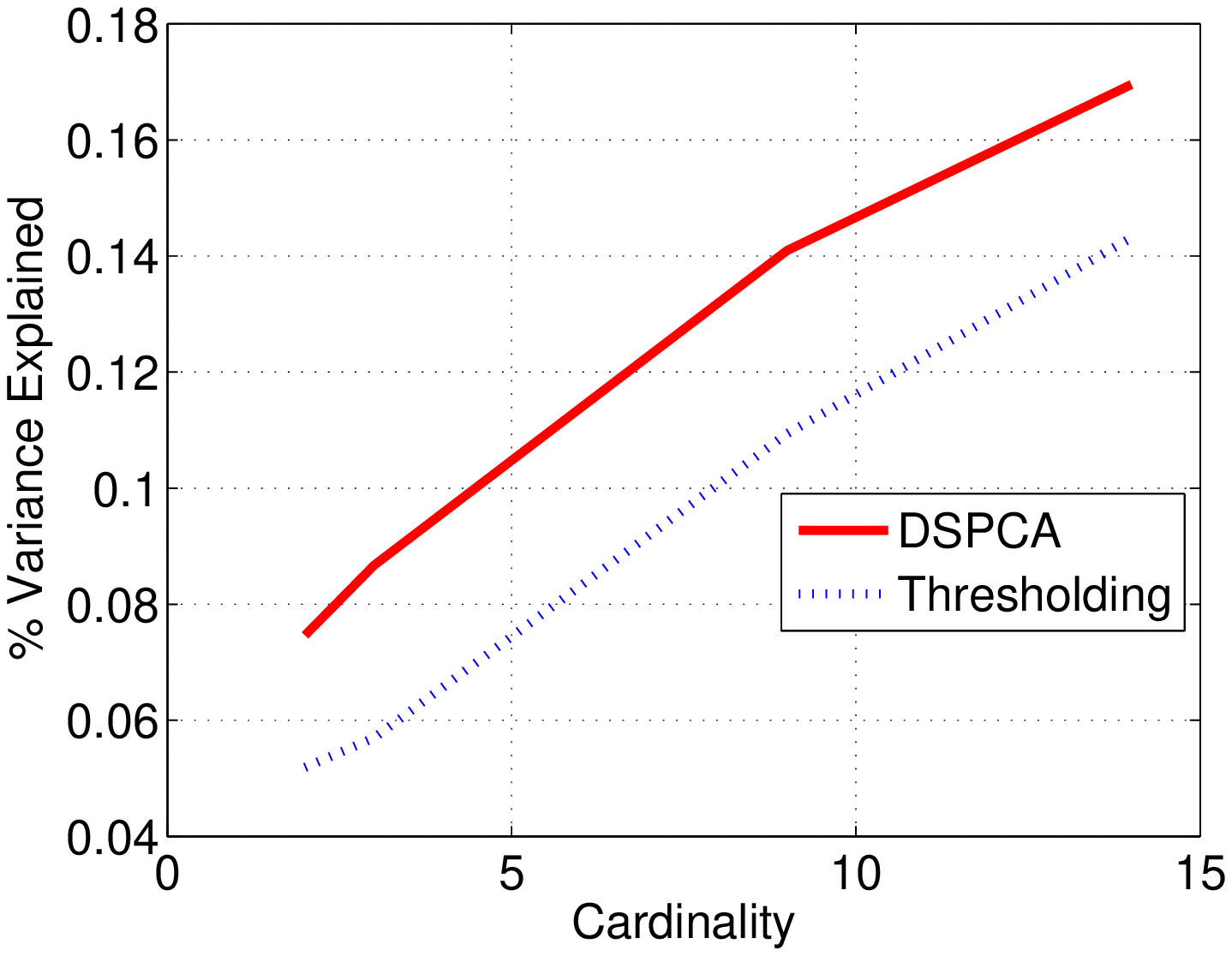}&
\psfrag{Cardinality}[t][b]{Cardinality}
\psfrag{\% Variance Explained}[b][t]{\% Variance Explained}
\includegraphics[width = 0.48\textwidth, height = 0.38\textwidth]{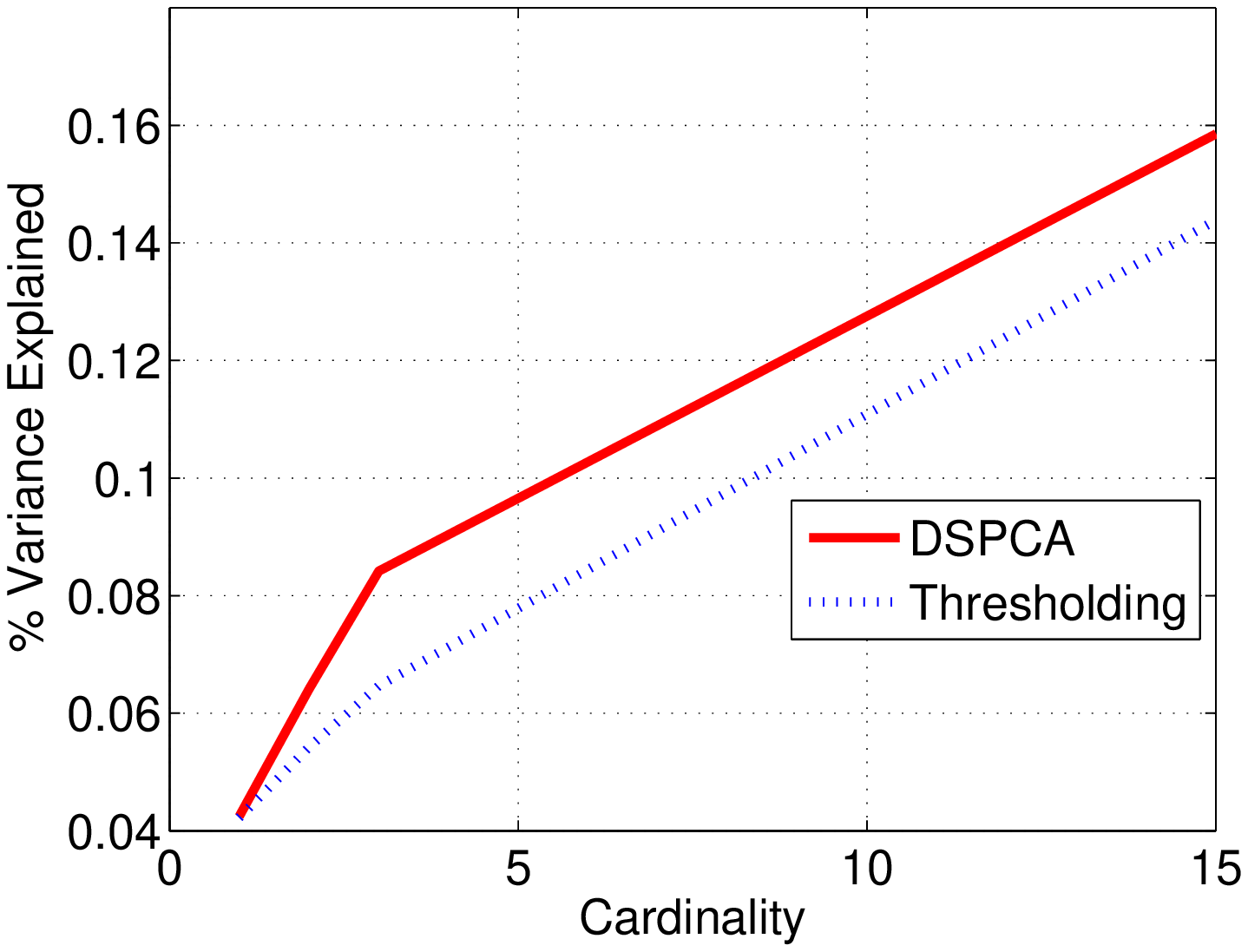}
\end{tabular}
\end{center}
\caption{Sparse PCA on 1,288 {\em{New York Times}} articles mentioning the word ``China''. \label{nytpc}}
\end{figure}

We also run the thresholded PCA and DSPCA algorithms over the collection of 1,288 news articles published by the {\em{New York Times}}'s International section mentioning the word ``China.'' We tokenize the articles by unigrams, remove no stop words, and perform no stemming. The data encodes the binary $\{0,1\}$ values (corresponding to appearance/non-appearance) of 86500 tokens. 

Figure~\ref{nytpc} shows the percentage of explained variance as a function of cardinality. Here we see DSPCA does outperform Thresholded PCA, though not by a big margin. Although we do not have ground truth, Table \ref{spca1pcwords} and Table \ref{pca1pcwords} contains words selected by two algorithms respectively as we increase cardinality. Words selected by DSPCA appear much more meaningful than those chosen by thresholded PCA at the same cardinality. 

\begin{table}[t]
\begin{center}\small{
\begin{tabular}{cccc}
$k=2$ & $k=3$ & $k=9$ & $k=14$ \\
\hline
\texttt{united} & \texttt{~american~} & \texttt{washington}  & \texttt{international}  \\
\texttt{states} & \texttt{united}  & \texttt{american} & \texttt{would}\\
& \texttt{states} & \texttt{~administration~} & \texttt{will}\\
& & \texttt{united} & \texttt{washington}\\
&  & \texttt{states}  & \texttt{american}\\
 &  & \texttt{president}  & \texttt{administration}\\
 && \texttt{obama}  &\texttt{united} \\
 & &\texttt{countries}   & \texttt{states}\\
 & &\texttt{nations}   & \texttt{president}\\
& &  & \texttt{obama}\\
 & &  & \texttt{counties}\\
 & &  & \texttt{nations}\\
  & &  & \texttt{policy}\\
& &  & \texttt{nuclear}\\
\hline \\
\end{tabular}}
\caption{1st PC from DSPCA on 1,288 {\em{New York Times}} articles mentioning the word ``China'' for various values of the eigenvector cardinality $k$.\label{spca1pcwords}}
\end{center}
\end{table}

\begin{table}[t!]
\begin{center}\small{
\begin{tabular}{cccc}
$k=2$ & $k=3$ & $k=9$ & $k=14$ \\
\hline
\texttt{even} & \texttt{even} & \texttt{even}  & \texttt{would} \\
\texttt{like} & \texttt{like}  & \texttt{we} & \texttt{new}\\
& \texttt{~states~} & \texttt{like} & \texttt{even}\\
& & \texttt{now} & \texttt{we}\\
&  & \texttt{this}  & \texttt{like}\\
 &  & \texttt{will}  & \texttt{now}\\
 && \texttt{united}  &\texttt{this} \\
 & &\texttt{~states~}   & \texttt{will}\\
 & &\texttt{if}   & \texttt{united}\\
& &  & \texttt{~states~}\\
 & &  & \texttt{world}\\
 & &  & \texttt{so}\\
  & &  & \texttt{some}\\
& &  & \texttt{if}\\
\hline \\
\end{tabular}}
\caption{1st PC from Thresholded PCA on 1,288 {\em{New York Times}} articles mentioning the word ``China'' for various values of the eigenvector cardinality $k$.\label{pca1pcwords}}
\end{center}
\end{table}

\subsection{Senate voting data}
In this section, we analyze the voting records of the 109th US Senate (2000-2004) 
There were 101 senators (one extra Senator is due to a replacement during the term) and 48  bills involved. To simplify, the votes are divided into ÒyesÓ (coded as 1) or ÒnoÓ (coded as -1), and  other votes are coded as 0.  

\begin{figure}[htb]
\begin{center}
\begin{tabular}{cc}
\psfrag{second principal component (cardinality = 48 )}[b][t]{2nd P.C. (Card = 48)}
\psfrag{first principal component (cardinality = 48 )}[t][b]{1st P.C. (Card = 48)}
\includegraphics[width = 0.48\textwidth, height = 0.38\textwidth]{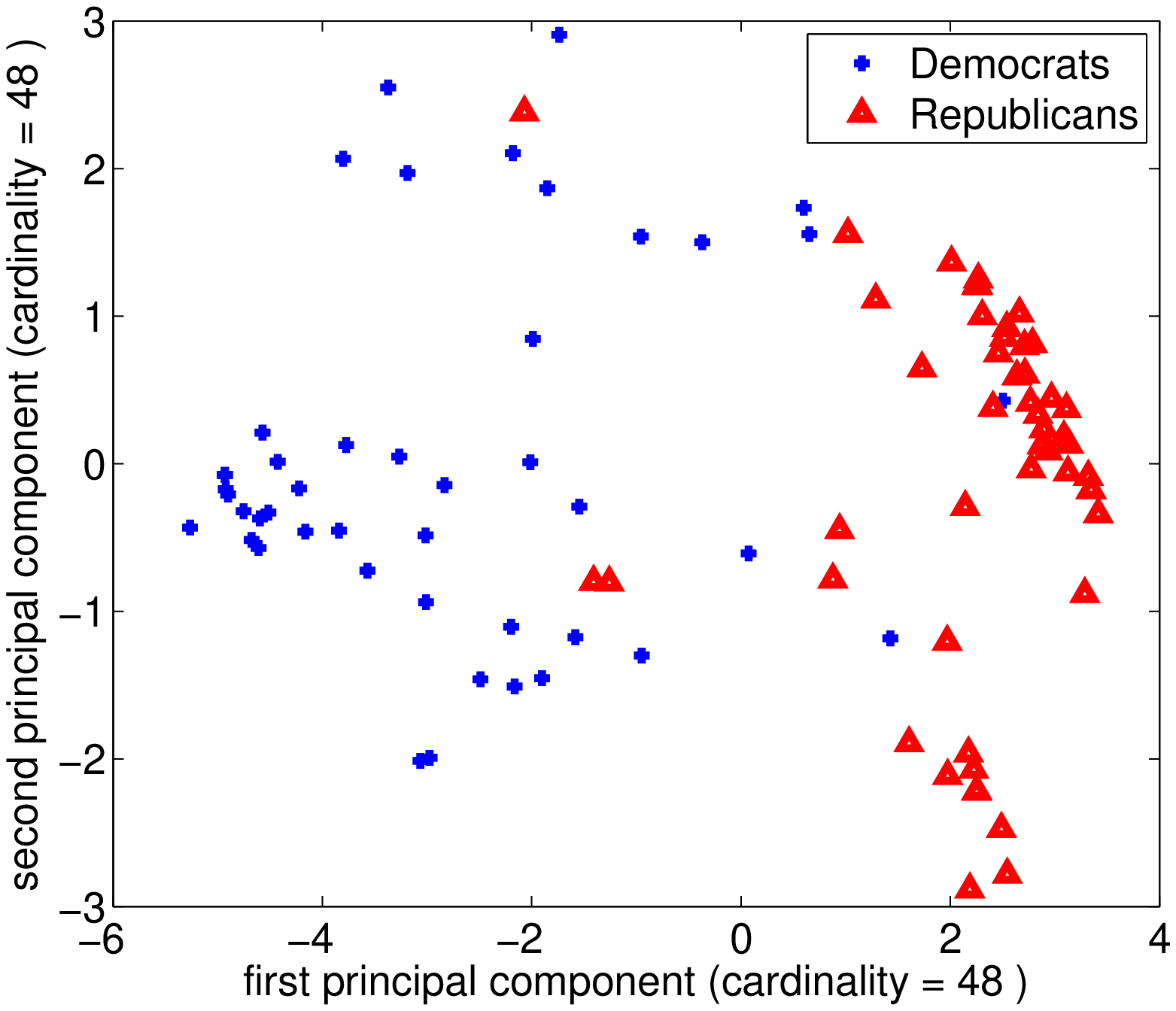}&
\psfrag{second principal component (cardinality = 2 )}[b][t]{2nd P.C. (Card = 2)}
\psfrag{first principal component (cardinality = 5 )}[t][b]{1st P.C. (Card = 5)}
\includegraphics[width = 0.48\textwidth, height = 0.38\textwidth]{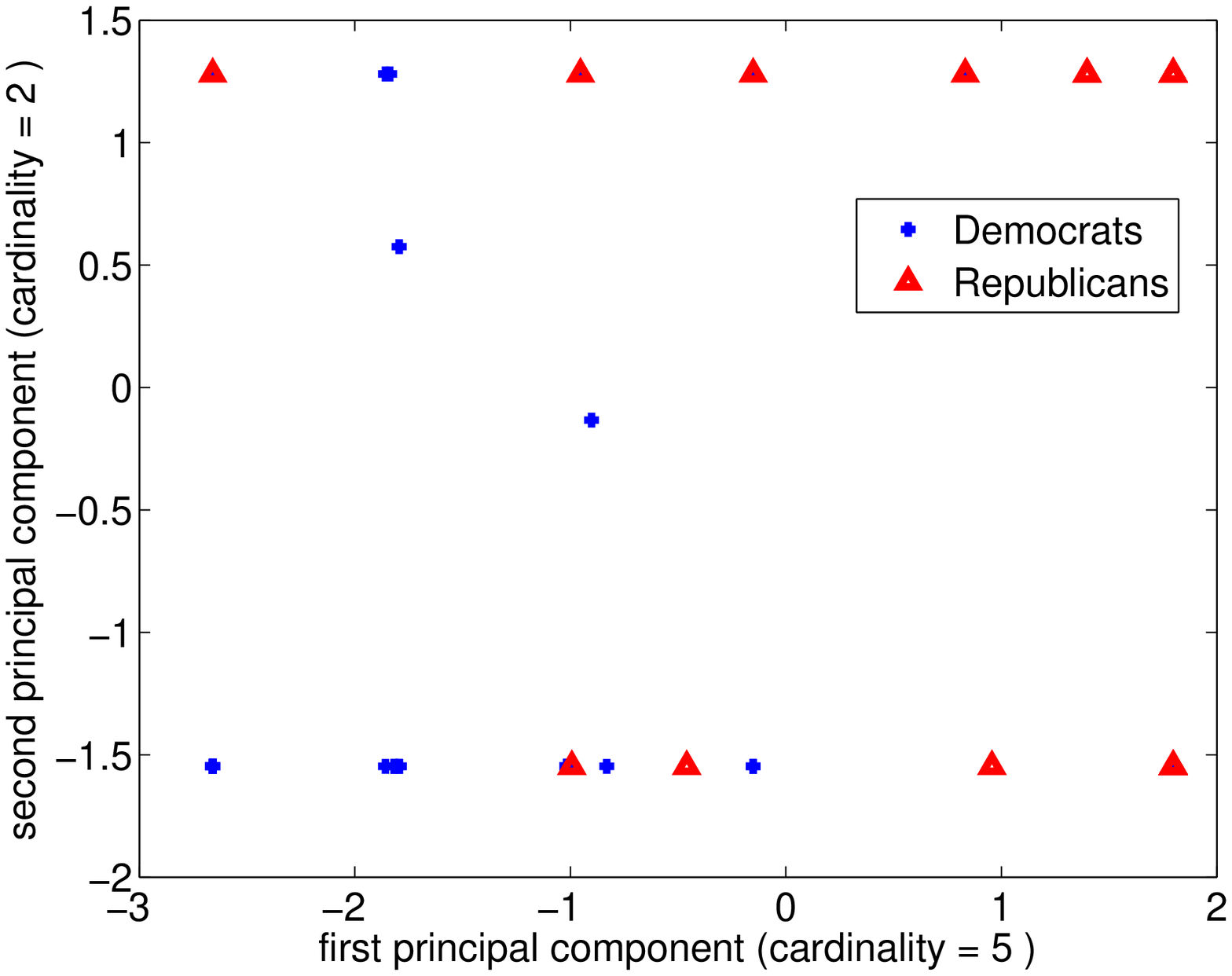}
\end{tabular}
\end{center}
\caption{109th Senate's voting record projected onto the top 2 principal components. \label{fig:svpca} }
\end{figure}

Each senator's voting record can be viewed as a point in a 48-dimensional space. By applying PCA, and projecting each senator's voting record onto a two-dimensional subspace of maximum variance, we can see that senators are almost perfectly separated by partisanship (Fig.~\ref{fig:svpca}). However, since the principal components involve all the bills, it is hard to tell which bills are most responsible for the explained variance. By applying Sparse PCA to the voting record, we aim to find a few bills that not only divide the senators according to partisanship, but also reveal which topics are most controversial within the Republican and Democratic parties. Fig.~\ref{fig:svpca} (right) shows the senators' voting records, projected onto the first two sparse principal components. We note that in the two-dimensional space senators are still divided by partisanship. In fact, many republican senators perfectly coincide with each other and so are democratic senators. In contrast to Fig.~\ref{fig:svpca} (left),  the cardinalities associated with the first and second sparse principal components are $5$ and $2$ respectively, which makes it possible to interpret the coordinates.

Let us examine the bills appearing in the first two sparse principal components. For the first sparse PC, the corresponding bills' brief description is as follows:
\begin{itemize}
\item S.\ 1932, As Amended; Deficit Reduction Act of 2005.

\item S.\ Con.\ Res.\ 83; An original concurrent resolution setting forth the congressional budget for the United States Government for fiscal year 2007 and including the appropriate budgetary levels for fiscal years 2006 and 2008 through 2011.

\item S.\ 3930, As Amended; Military Commissions Act of 2006.

\item S.\ 403, As Amended; Child Interstate Abortion Notification Act.

\item Passage of S.\ 397, As Amended; Protection of Lawful Commerce in Arms Act.
\end{itemize}

The brief description for the two bills in the second sparse principal component are:
\begin{itemize}
\item H.\ R.\ 3045; Dominican Republic-Central America-United States Free Trade Agreement Implementation Act.

\item S.\ 1307; Dominican Republic-Central America-United States Free Trade Agreement Implementation Act.
\end{itemize}
A glance at these bills tells us that the major controversial issues between Democrats and Republicans are topics such as ``abortion", ``military", ``budget", and ``free trade". 

\begin{figure}[htb]
\begin{center}
\begin{tabular}{cc}
\psfrag{Cardinality}[t][b]{Cardinality}
\psfrag{Variance}[b][t]{\% Variance Explained}
\includegraphics[width = 0.48\textwidth, height = 0.38\textwidth]{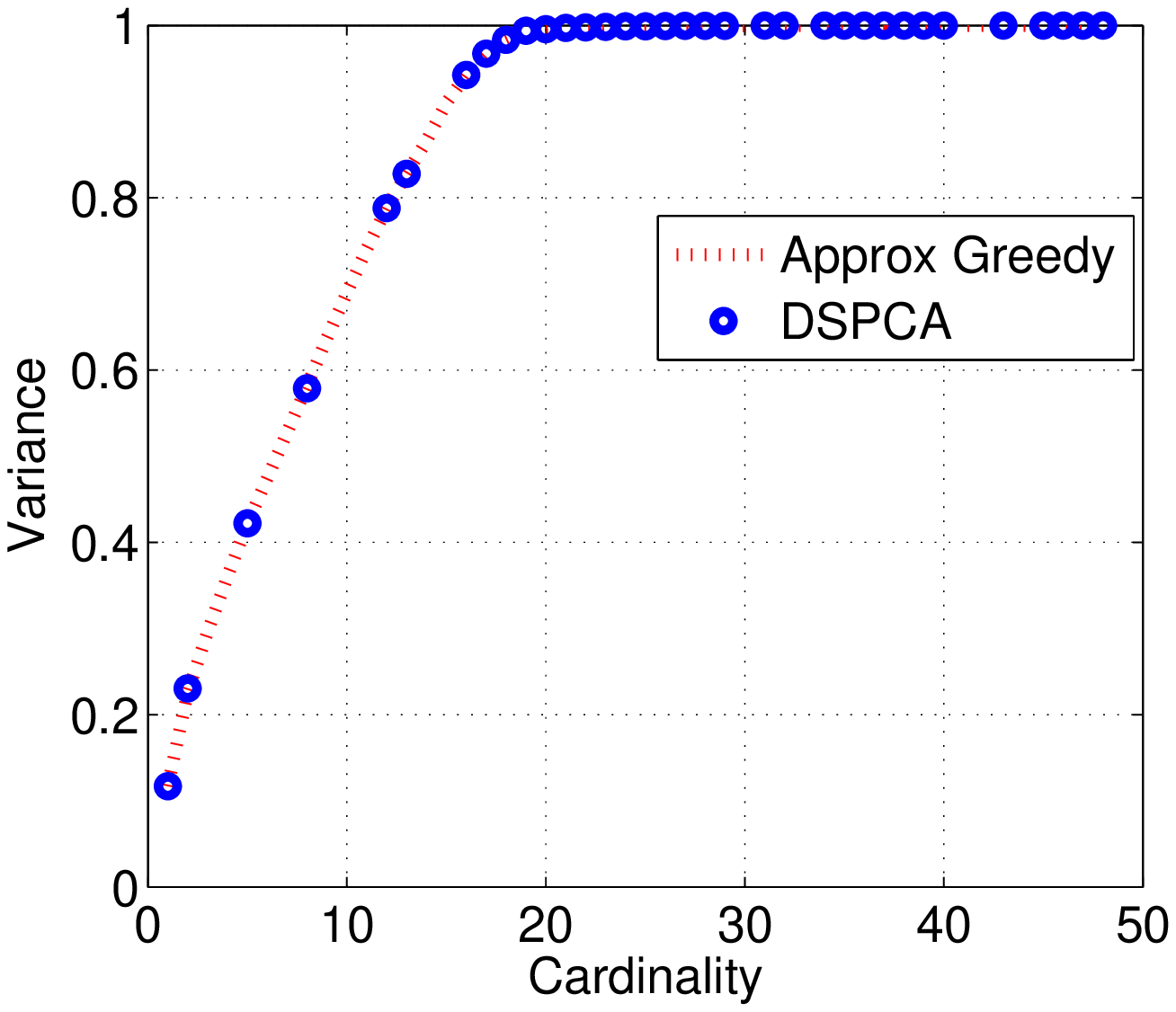}&
\psfrag{\rho}[b][t]{\% Variance Explained}
\psfrag{Cardinality}[b][t]{Cardinality}
\includegraphics[width = 0.48\textwidth, height = 0.38\textwidth]{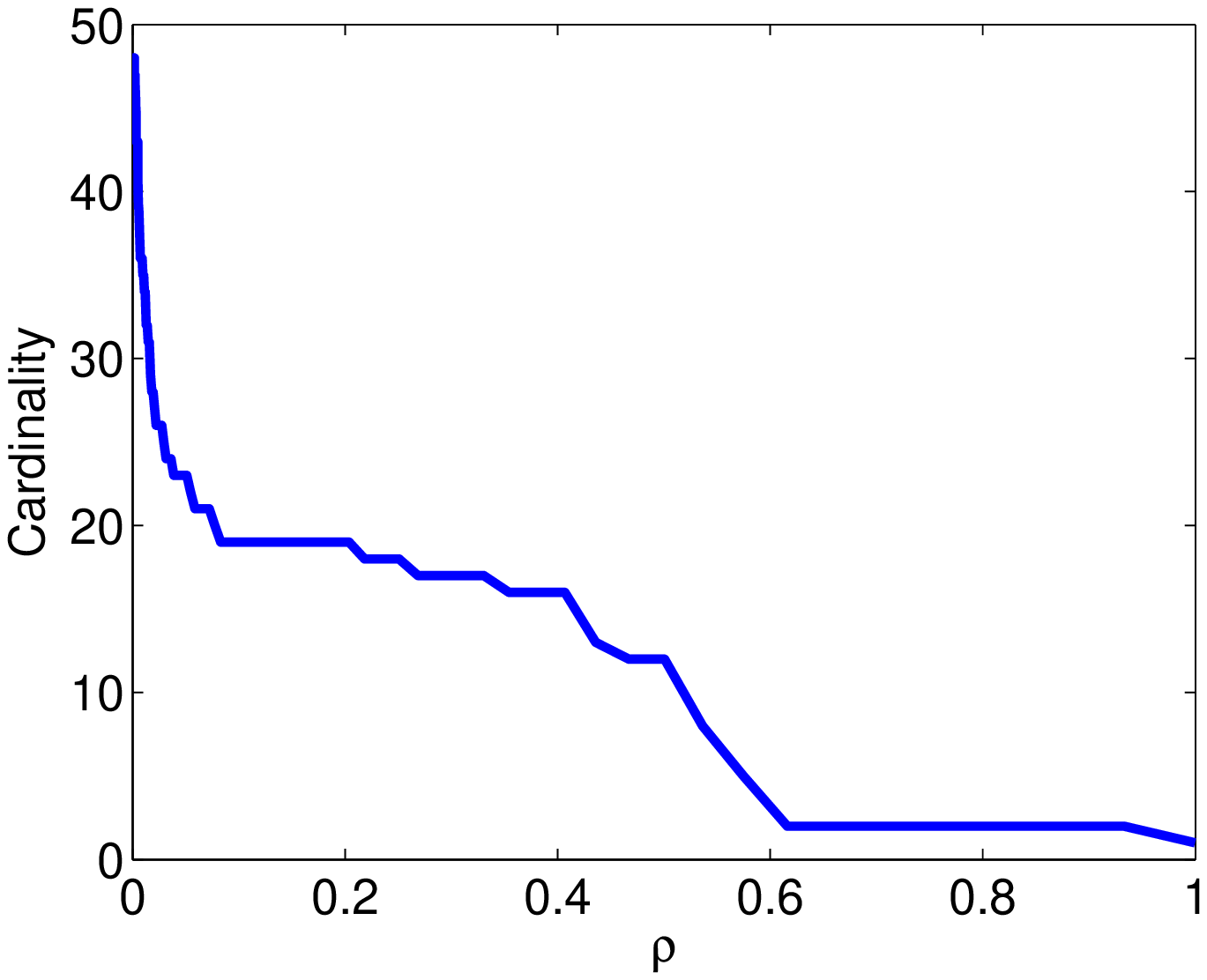}
\end{tabular}
\caption{{\em Left:} Explained variance as a function of cardinality. {\em Right:} Cardinality as a function of penalty parameter $\rho$. \label{fig:svexpvar}}
\end{center}
\end{figure}

Fig ~\ref{fig:svexpvar} plots the variance explained by the first sparse principal component divided by that explained by the first PC, as a function of the cardinality of the sparse PC. Fig.~\ref{fig:svexpvar} also shows how the cardinality of the first sparse PC varies as the penalty parameter $\rho$ is changed in the DSPCA code.We can see that when 19 out of 48 variables (bills) are used, sparse PCA almost achieves the same statistical fidelity as standard PCA does. 
 
\subsection{Stock market data}
In this section, we investigate the historical prices of S\&P500 stocks over 5 years, from June 1st, 2005, through June 1st, 2010. By taking out the stocks with less than 5 years of history, we end up with 472 stocks, each having daily closing prices over 1259 trading days. The prices are first adjusted for dividends and splits and then used to calculate daily log returns. Each day's return can be represented as a point in $\reals^{472}$. 

\begin{figure}[htb]
\begin{center}
\begin{tabular}{cc}
\psfrag{Cardinality}[t][b]{Cardinality}
\psfrag{Variance}[b][t]{\% Variance Explained}
\includegraphics[width = 0.48\textwidth, height = 0.38\textwidth]{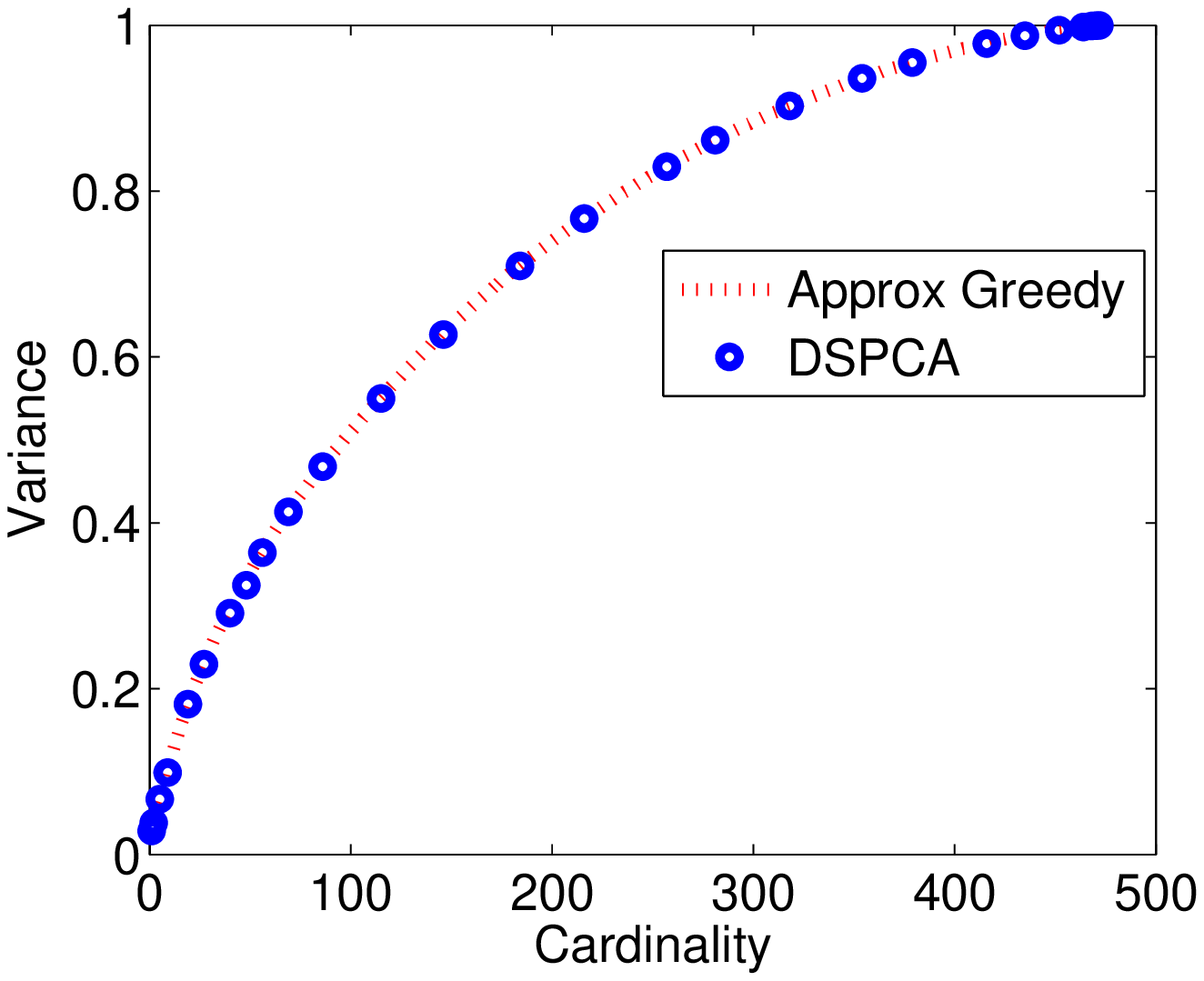}&
\psfrag{\rho}[b][t]{\% Variance Explained}
\psfrag{Cardinality}[b][t]{Cardinality}
\includegraphics[width = 0.48\textwidth, height = 0.38\textwidth]{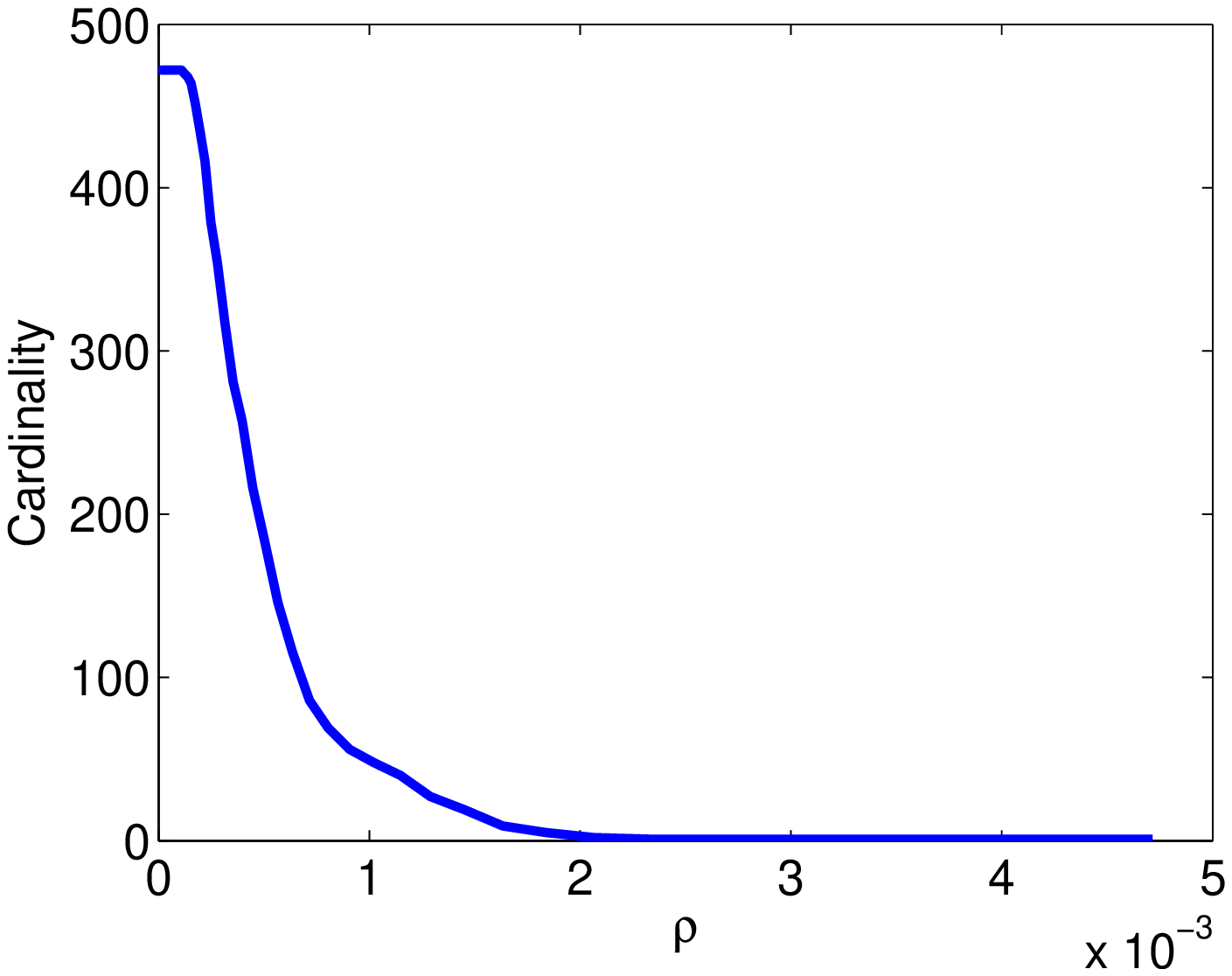}
\end{tabular}
\caption{{\em Left:} Explained variance as a function of cardinality. {\em Right:} Cardinality as a function of penalty parameter $\rho$. \label{fig:spexpvar}}
\end{center}
\end{figure}

Fig.~\ref{fig:spexpvar} shows the explained variance as a function of 1st PC's cardinality. It seems hard to say that the 1st PC is sparse, since there is no natural ``kink" in that curve.  That is, we need almost 300 out of the total 472 stocks to explain at least $90\%$ of the variance explained by the 1st PC from PCA. However, when we inspect the sparse PCs with increasing cardinalities, we note that initially only stocks from the "Financials" sector come to play and later until, at cardinality 32, do we see companies from other sectors appearing in the 1st sparse PC. So we take the first sparse PC with cardinality equal to 32. Then we solve for the 2nd sparse PC, and using the same guideline to arrive at a cardinality of 26. 

\begin{figure}[htb]
\begin{center}
\begin{tabular}{cc}
\psfrag{second principal component (cardinality = 472 )}[b][t]{2nd P.C. (Card = 472)}
\psfrag{first principal component (cardinality = 472 )}[t][b]{1st P.C. (Card = 472)}
\includegraphics[width = 0.48\textwidth, height = 0.38\textwidth]{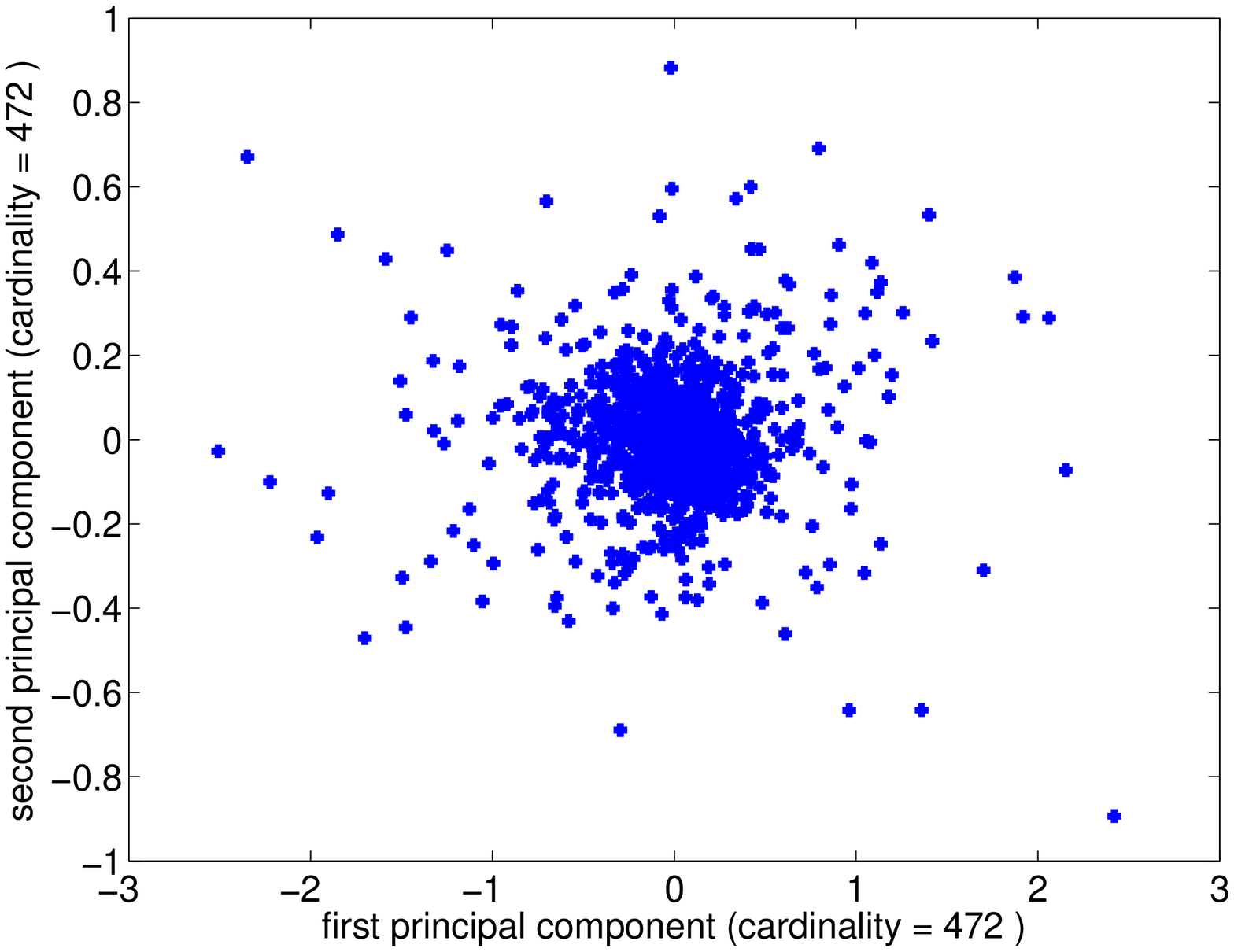}&
\psfrag{second principal component (cardinality = 26 )}[b][t]{2nd P.C. (Card = 26)}
\psfrag{first principal component (cardinality = 32 )}[t][b]{1st P.C. (Card = 32)}
\includegraphics[width = 0.48\textwidth, height = 0.38\textwidth]{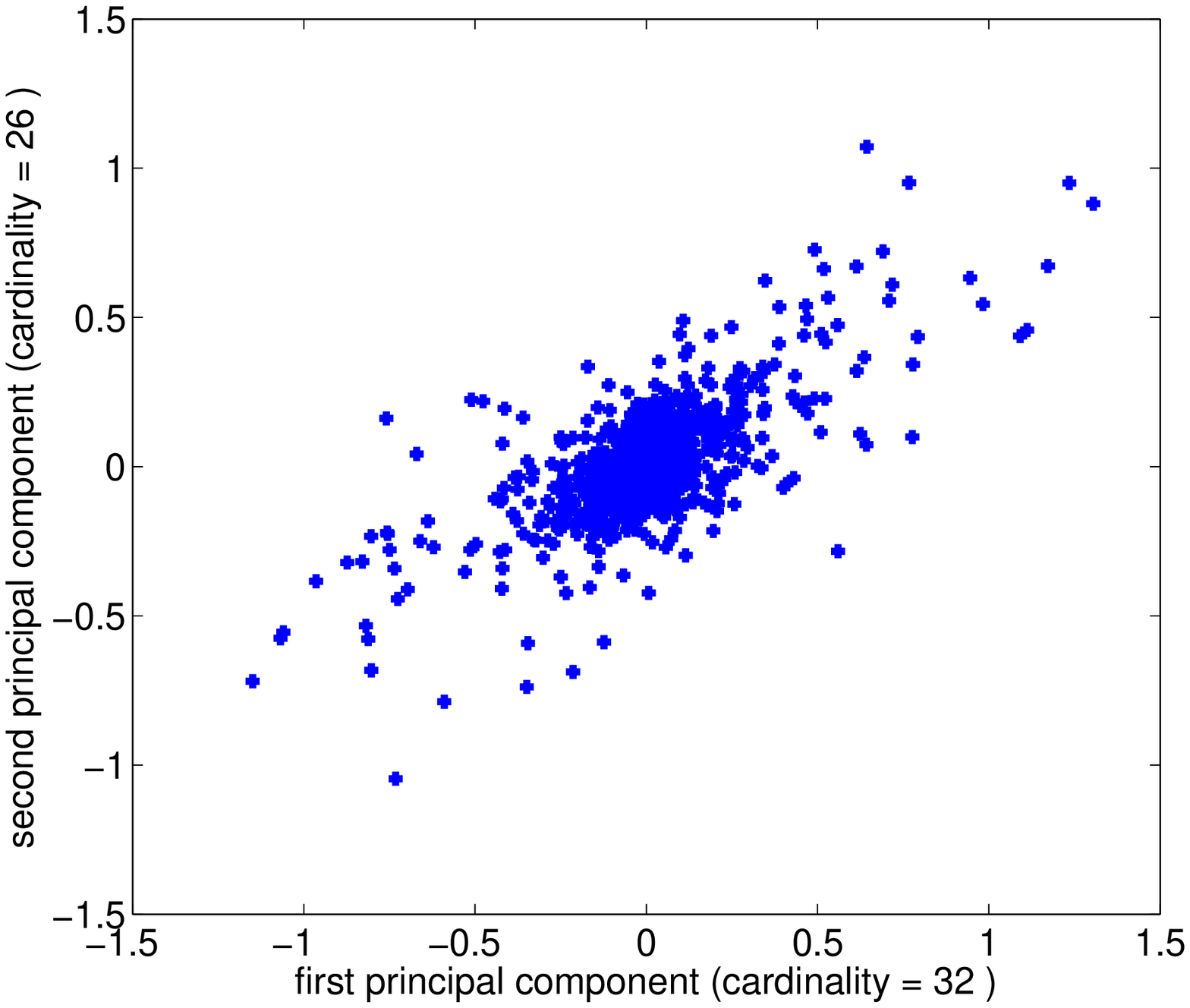}
\end{tabular}
\end{center}
\caption{S\&P500 daily returns projected onto the top 2 principal components. For PCA (left) and sparse PCA (right).\label{fig:sppca}}
\end{figure}

Figure~\ref{fig:sppca} show the stock returns projected onto the 2-dimensional subspaces spanned by the top 2 PCs and top 2 sparse PCs, respectively. Comparing these two plots, we observe two interesting phenomena:
\begin{itemize}
\item Although the top 2 principal components from PCA explain more variance (as seen from the larger range of the axes in the left over the right panel), the two sparse principal components from DSPCA involve only 58 out of 472 stocks (32 on the first PC and another distinct 26 on the second). Furthermore,  31 of the 32 stocks in the first sparse PC are all from the sector "Financials", and that almost all 26 stocks in the second sparse PC come from "Energy" and "Materials" except 2 from "Financials" and 1 from "Information Technology". Considering that there are 10 sectors in total, this is quite interesting as Sparse PCA is able to identify the right groups (industry factors) that explains most of the variance. Our data covers June 2005 through June 2010 where a severe financial crisis took place, and the key role of the Financial sector is revealed purely through our sparse PCA analysis.

\item In Fig.~\ref{fig:spexpvar}, the projected data appears symmetrically distributed around its center.  In contrast, In Fig.~\ref{fig:sppca} (right),  we observe a definite orientation.  Since the horizontal axis (first PC) corresponds to ``Financials" and the vertical one to ``Energy" and ``Materials", the sparse PCA analysis tells us that these two sectors are positively correlated.
\end{itemize}

\subsection{Random matrices}
\begin{figure}
\begin{center}
\psfrag{card}[t][b]{Cardinality}
\psfrag{var}[b][t]{Max. Eigenvalue}
\includegraphics[width=.48\textwidth]{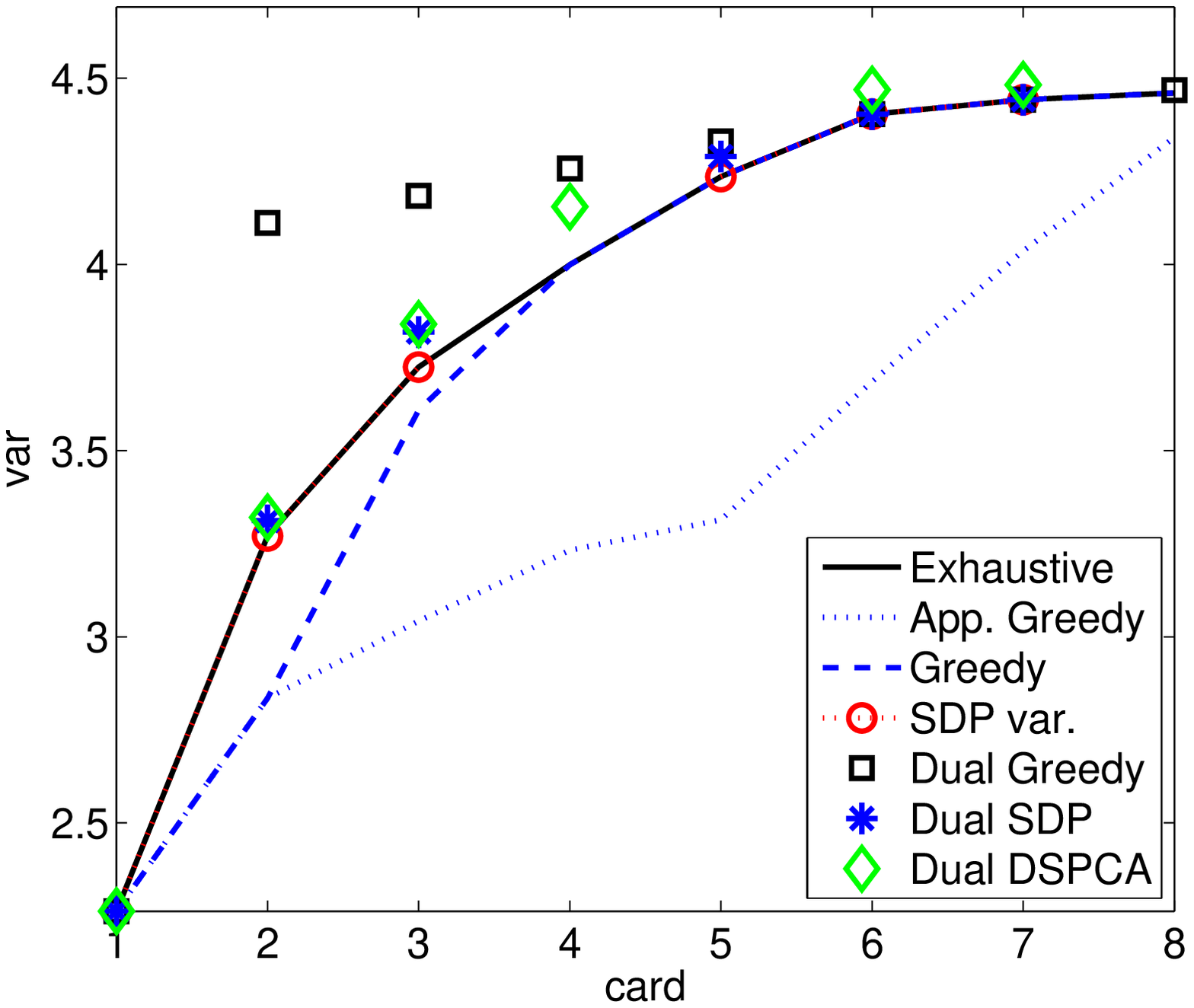}
\includegraphics[width=.48\textwidth]{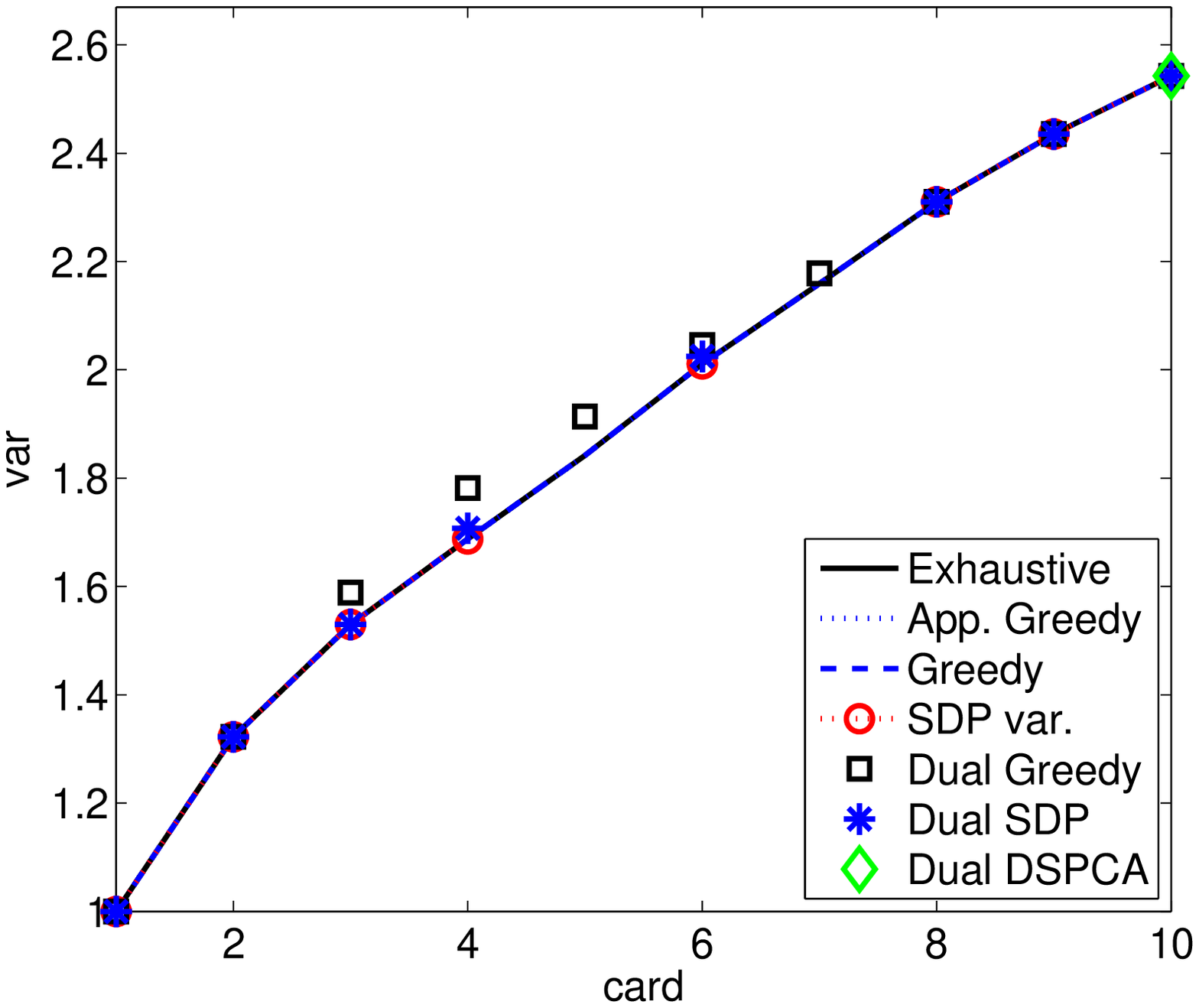}
\end{center}
\caption{Upper and lower bound on sparse maximum eigenvalues. We plot the maximum sparse eigenvalue versus cardinality, obtained using exhaustive search (solid line), the approximate greedy (dotted line) and fully greedy (dashed line) algorithms. We also plot the upper bounds obtained by minimizing the gap of a rank one solution (squares), by solving the $\ell_0$ semidefinite relaxation explicitly (stars) and by solving the DSPCA dual $\ell_1$ relaxation (diamonds). \emph{Left:} On a matrix $F^TF$ with $F$ Gaussian. \emph{Right:} On a sparse rank one plus noise matrix.\label{fig:recov}}
\end{figure}

Sparse eigenvalues of random matrices play a central role in characterizing the performance of $\ell_1$ decoders in compressed sensing applications. Testing the Restricted Isometry Property (RIP) in \cite{Cand05} amounts to bounding the maximum and minimum eigenvalues of a Gram matrix. Here, we compute the upper and lower bounds on sparse eigenvalues produced using various algorithms. We pick the data matrix to be small enough so that computing sparse eigenvalues by exhaustive search is numerically feasible. In Figure~\ref{fig:recov}, we plot the maximum sparse eigenvalue versus cardinality, obtained using exhaustive search (solid line), the approximate greedy (dotted line) and fully greedy (dashed line) algorithms. We also plot the upper bounds obtained by minimizing the gap of a rank one solution (squares), by solving the semidefinite relaxation in~\S\ref{ss:l0} explicitly (stars) and by solving the DSPCA dual (diamonds). On the left, we use a matrix $\Sigma=F^TF$ with $F$ Gaussian. On the right, $\Sigma=uu^T/\|u\|^2+2V^TV$, where $u_i=1/i,~i=1,\ldots,n$ and $V$ is matrix with coefficients uniformly distributed in $[0,1]$. Almost all algorithms are provably optimal in the noisy rank one case (as well as in many example arising from ``natural data''), while Gaussian random matrices are harder. Note however, that the duality gap between the semidefinite relaxations and the optimal solution is very small in both cases, while our bounds based on greedy solutions are not as good. Overall, while all algorithms seem to behave similarly on ``natural'' or easy data sets, only numerically expensive relaxations produce good bounds on the random matrices used in compressed sensing applications.

\subsection{Statistical consistency vs. computational complexity}\label{ss:consitency}
As we hinted above, very simple methods such as thresholding or greedy algorithms often perform well enough on simple data sets, while obtaining good statistical fidelity on more complex (or random) data sets requires more complex algorithms. This is perfectly illustrated by the results in \cite{Amin08} on a spiked covariance model. To summarize these results, suppose that the sample covariance matrix $\hat \Sigma \in \symm_n$ is a noisy estimate of the true population covariance $\Sigma\in\symm_n$ with $\hat \Sigma = \Sigma + \Delta$ where $\Delta$ is a noise matrix, suppose also that the leading eigenvector of the true covariance is sparse with cardinality $k$. Under some assumptions on the noise component $\Delta$, \cite{Amin08} show that when the ambient dimension $n$, the number of observations $m$ and the number $k$ of nonzero components in the leading eigenvector all scale to infinity, and when the ratio
\[
\theta_\mathrm{thres}=\frac{m}{k^2 \log(n-k)}
\]
is above some critical value, then simply thresholding the diagonal of the sample covariance matrix will recover the exact support of the leading eigenvector of $\Sigma$ with probability tending to one. On the other hand, simple thresholding fails with probability one when this ratio is below a certain value. Furthermore, \cite{Amin08} show that when
\[
\theta_\mathrm{sdp}=\frac{m}{k \log(n-k)}
\]
is above some critical value, the solution of the semidefinite relaxation in Section~\ref{ss:l1-relax} (if tight) will recover the exact support of the leading eigenvector of $\Sigma$ with probability tending to one. On the other hand, the semidefinite relaxation fails with probability one when this ratio is below a certain value. They also show that the semidefinite programing relaxation in Section~\ref{ss:l1-relax} is statistically optimal, meaning that no other method (even combinatorial ones) can recover the true support using fewer samples (up to a constant factor). This result clearly illustrates a tradeoff between statistical fidelity on one side and computational complexity on the other. In the spiked model, the semidefinite relaxation requires $O(1/k)$ fewer samples than simply thresholding the diagonal to recover the true support of the leading eigenvector of $\Sigma$, but its complexity is much higher than that of the thresholding strategy.

We can further illustrate this behavior on a simple numerical example. Suppose we are given a sample covariance $\hat \Sigma \in \symm_n$ coming from a ``spiked'' model of covariance similar to that in \cite{Amin08}, with
\[
\hat \Sigma = uu^T + VV^T/\sqrt{m}
\]
where $u\in\reals^n$ is the true sparse leading eigenvector, with $\Card(u)=k$, $V\in\reals^{n \times m}$ is a noise matrix with $V_{ij}\sim{\cal N}(0,1)$ and $m$ is the number of observations. We compare the performance of the simple thresholding method (on the leading eigenvector of regular PCA here) with that of the semidefinite relaxation when recovering the support of $u$ for various values of the number of samples. Our point here is that, while variance versus cardinality is a direct way of comparing the  performance of sparse PCA algorithms, accurate recovery of the support is often a far more important objective. Many methods produce similar variance levels given a limited budget of nonzero components, but their performance in recovering the true support is often markedly different. 

\begin{figure}[htb]
\begin{center}
\begin{tabular}{cc}
\psfrag{spec}[t][b]{Specificity}
\psfrag{sens}[b][t]{Sensitivity}
\includegraphics[width = 0.48\textwidth, height = 0.38\textwidth]{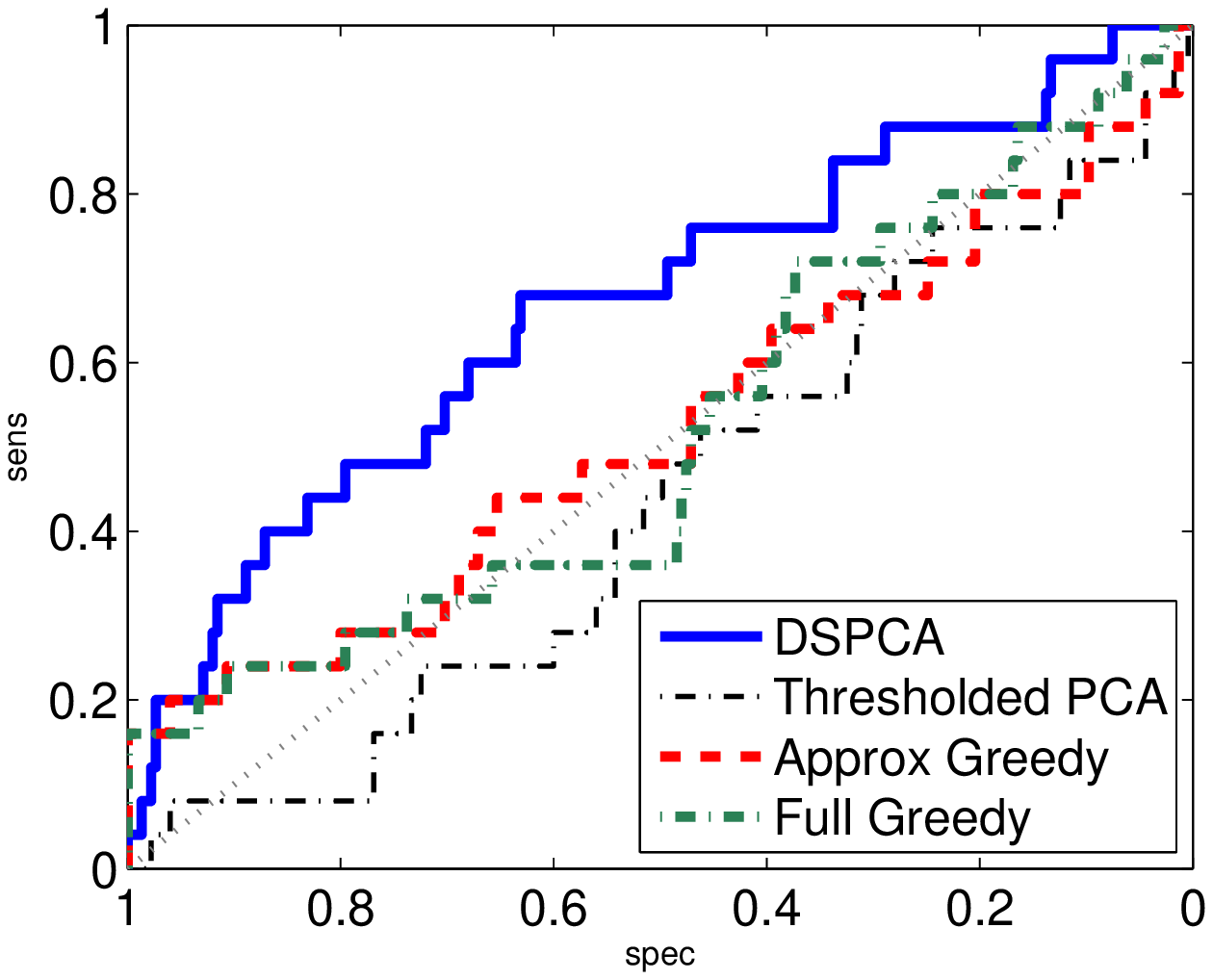}&
\psfrag{Number of Samples}[t][b]{Number of samples $m$}
\psfrag{AUROC}[b][t]{Area Under ROC}
\includegraphics[width = 0.48\textwidth, height = 0.38\textwidth]{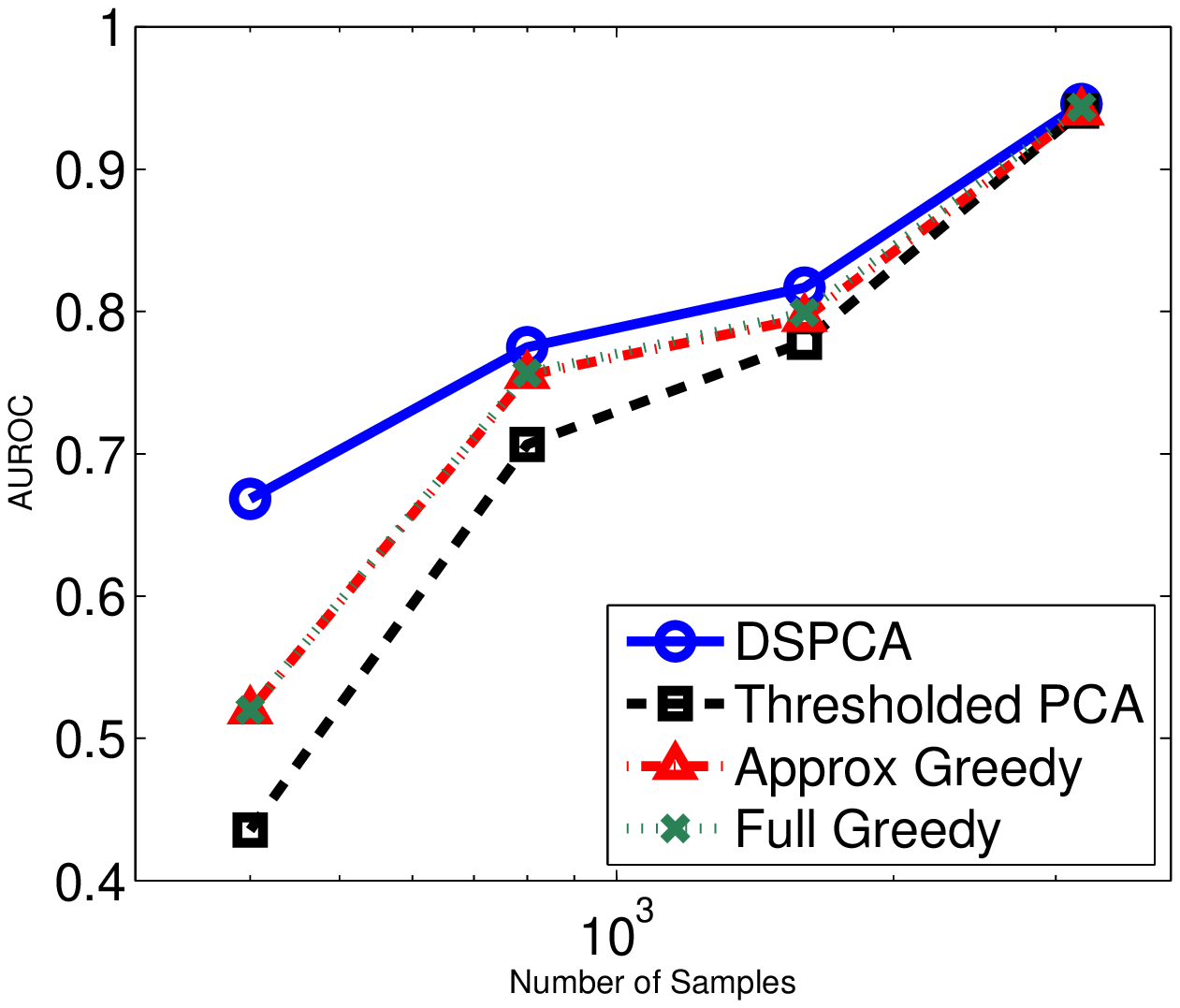}
\end{tabular}
\caption{{\em Left:} ROC curves when recovering the support of $u$ in the spiked model using thresholding, approximate and exact greedy algorithms and the semidefinite relaxation (DSPCA) in Section~\ref{ss:l1-relax} in the spiked model when $n=250$, $m=400$ and $k=25$. {\em Right:}  Area Under ROC (AUROC) versus number of samples $m$.\label{fig:roc}}
\end{center}
\end{figure}

In Figure~\ref{fig:roc} on the left we compare ROC curves when recovering the support of $u$ in the spiked model above using thresholded PCA, the approximate and full greedy algorithms in \cite{dAsp08b} and semidefinite relaxation (DSPCA). On the right, we plot Area Under ROC as the number of samples increase. As expected, we observe that the semidefinite relaxation performs much better when only a limited number of observations are available ($m$ small). 

\section{Conclusion}
We have reviewed here several techniques for approximating the solution to the {\em single factor} sparse PCA problem. While the algorithms presented here perform quite well, several key questions remain open at this point.

First, outside of the (locally) convergent algorithm in \cite{Jour08}, very few methods handle the problem of simultaneously finding several leading sparse principal components.

Also, as the examples of Section~\ref{s:apps} illustrate, most methods (even extremely simple ones) perform well enough on easy, ``natural'' data sets while only the most expensive semidefinite relaxations seem to produce good bounds on the random matrices used in compressed sensing applications for example. Characterizing what makes ``natural'' data sets easier than random ones remains an open problem at this point. It is also not clear yet how to extend the statistical optimality statements of \cite{Amin08} to broader (e.g. deterministic) classes of matrices.

Finally, the question of approximation bounds \`a la MAXCUT for the relaxations detailed here is largely open. Basic performance bounds are discussed in \cite{dAsp08a, Bach10} but they can certainly be tightened.

\section*{Acknowledgments} The authors gratefully acknowledge partial support from NSF grants SES-0835550 (CDI), CMMI-0844795 (CAREER), CMMI-0968842, a Peek junior faculty fellowship, a Howard B. Wentz Jr. award and a gift from Google. 

\small{\bibliographystyle{plainnat}
\bibsep 1ex
\bibliography{MainPerso}}
\end{document}